\documentclass[12pt]{amsart}
\usepackage{amssymb}
\usepackage{amscd}
\usepackage{thmdefs}

\input{tcilatex}
\theoremstyle{definition}
\theoremstyle{remark}
\numberwithin{equation}{section}

\input tcilatex

\begin{document}
\title[A Curvature Flow on a CR $3$-manifold]{A Fourth Order Curvature Flow
on a CR $3$-manifold }
\author{$^{\ast }$Shu-Cheng Chang$^{1}$}
\address{$^{1}$Department of Mathematics, National Tsing Hua University,
Hsinchu 30013, Taiwan, R.O.C.}
\email{scchang@math.nthu.edu.tw }
\thanks{$^{\ast }$Research supported in part by the NSC of Taiwan}
\author{$^{\ast }$Jih-Hsin Cheng$^{2}$}
\address{$^{2}$Institute of Mathematics, Academia Sinica, Nankang, Taipei
11529, Taiwan, R.O.C.}
\email{cheng@math.sinica.edu.tw}
\author{Hung-Lin Chiu$^{3}$}
\address{$^{3}$Institute of Mathematics, Academia Sinica, Nankang, Taipei
11529, Taiwan, R.O.C.}
\email{hlchiu@math.sinica.edu.tw}

\begin{abstract}
Let $(\mathbf{M}^{3},J,\theta _{0})$ be a closed pseudohermitian $3$-manifold%
$.$ Suppose the associated torsion vanishes and the associated $Q$-curvature
has no kernel part with respect to the associated Paneitz operator. On such
a background pseudohermitian $3$-manifold, we study the change of the
contact form according to a certain version of normalized $Q$-curvature
flow. This is a fourth order evolution equation. We prove that the solution
exists for all time and converges smoothly to a contact form of zero $Q$%
-curvature. We also consider other background conditions and obtain a priori
bounds up to high orders for the solution.
\end{abstract}

\maketitle

\section{\textbf{Introduction}}

Let $(M$, $J,$ $\theta )$ be a closed (i.e., compact with no boundary)
pseudohermitian $3$-manifold (see Appendix A for basic notions in
pseudohermitian geometry)$.$ In the papers \cite{fh}, \cite{gg}, \cite{h},
and \cite{gl}, the $(CR)$ Paneitz operator $P$ (acting on a smooth real
function $\lambda )$ with respect to $(J,$ $\theta )$ is defined by%
\begin{equation*}
P\lambda =\Delta _{b}^{2}\lambda +T^{2}\lambda +4\func{Im}(A_{\overline{11}%
}\lambda _{11}+A_{\overline{11},1}\lambda _{1}),
\end{equation*}%
and the so-called $Q$-curvature is defined by 
\begin{equation}
Q=-\frac{2}{3}(\Delta _{b}W+2\func{Im}A_{11\overline{,11}})  \label{-1}
\end{equation}%
where $\Delta _{b},$ $T,$ $W,\ $and $A_{\overline{11}}$ denote the
sub-Laplacian, the characteristic vector field, the Tanaka-Webster (scalar)
curvature, and the torsion with respect to $($ $J,$ $\theta ),$
respectively. Moreover, for a contact form change $\theta =e^{2\lambda
}\theta _{0},$ we have the following transformation laws:%
\begin{equation}
P=e^{-4\lambda }P_{0}  \label{0a}
\end{equation}%
and 
\begin{equation}
Q=e^{-4\lambda }(Q_{0}+2P_{0}\lambda )  \label{0b}
\end{equation}%
where $P_{0}$ and $Q_{0}$ denote the $(CR)$ Paneitz operator and the $Q$%
-curvature with respect to $(J,$ $\theta _{0}),$ respectively. Note that 
\begin{equation}
\int_{M}Qd\mu =0.  \label{1b}
\end{equation}%
Here the volume form $d\mu $ $=$ $\theta \wedge d\theta .$

Then we ask if we can always choose $\theta $ so that $Q$ vanishes
pointwise. The related problems for the $Q$-curvature on a Riemannian
manifold are also addressed and studied in \cite{b}, \cite{c1}, and \cite{cw}
by using the method of the $Q$-curvature flow. In this paper we study a
similar evolution equation and apply it to the problem addressed above on a $%
CR$ $3$-manifold.

We consider the functional $\mathcal{E}$ on a closed $CR$ $3$-manifold $%
(M,J) $ with a given contact class $[\theta _{0}]$ (which consists of all
contact forms annihilating the underlying contact bundle)$:$

\begin{equation}
\mathcal{E}\emph{(}\theta \emph{)}=\int_{M}P_{0}\lambda \cdot \lambda d\mu
_{0}+\int_{M}Q_{0}\lambda d\mu _{0},  \label{1a}
\end{equation}%
for $\theta =e^{2\lambda }\theta _{0},$ where $d\mu _{0}$ $=$ $\theta
_{0}\wedge d\theta _{0}.$ Then, for minimizing $\mathcal{E}(\theta )$ in $%
[\theta _{0}],$ it is natural to consider the following fourth order ($Q$%
-curvature) flow on a closed $CR$ $3$-manifold $(M,J):$%
\begin{equation}
\left\{ 
\begin{array}{l}
\frac{\partial \lambda }{\partial t}=-(Q_{0}+2P_{0}\lambda )+r=-e^{4\lambda
}Q+r, \\ 
\theta =e^{2\lambda }\theta _{0};\lambda (p,0)=\lambda _{0}(p), \\ 
\int_{M}e^{4\lambda _{0}}d\mu _{0}=\int_{M}d\mu _{0},%
\end{array}%
\right.  \label{2}
\end{equation}%
where%
\begin{equation*}
r=\frac{\int_{M}(Q_{0}+2P_{0}\lambda )d\mu }{\int_{M}d\mu }=\frac{%
\int_{M}e^{4\lambda }Qd\mu }{\int_{M}d\mu }
\end{equation*}%
and $\lambda _{0}$ is an initial real $C^{\infty }$ smooth function$.$ Note
that $d\mu =e^{4\lambda }d\mu _{0}$ and the volume $V=\int_{M}d\mu $ is kept
invariant under the flow (\ref{2}). Indeed (\ref{2}) is the (volume
normalized) negative gradient flow of $\mathcal{E}\emph{(}\theta \emph{).}$
That is, by using (\ref{1b}) and (\ref{0b}), one can check that 
\begin{equation}
\begin{array}{l}
\frac{d}{dt}\mathcal{E}\emph{(}\theta \emph{)} \\ 
=2\int_{M}P_{0}\lambda \cdot \frac{\partial \lambda }{\partial t}d\mu
_{0}+\int_{M}Q_{0}\frac{\partial \lambda }{\partial t}d\mu _{0} \\ 
=-2\int_{M}(Q_{0}+2P_{0}\lambda )P_{0}\lambda d\mu
_{0}-\int_{M}Q_{0}(Q_{0}+2P_{0}\lambda )d\mu _{0} \\ 
=-\int_{M}(Q_{0}+2P_{0}\lambda )^{2}d\mu _{0} \\ 
=-\int_{M}e^{4\lambda }Q^{2}d\mu .%
\end{array}
\label{2d}
\end{equation}%
under the flow (\ref{2}). We will often denote a $CR$ $3$-manifold by $%
(M,J,[\theta _{0}])$ with its contact class (or equivalently its underlying
contact bundle) indicated. Let $\overset{0}{A}_{11}$ denote the torsion with
respect to $(J,\theta _{0}).$ For the long time solution to the flow (\ref{2}%
), we have the following result.

\begin{theorem}
\label{Th1.1}Let $(M,J,[\theta _{0}])$ be a closed CR $3$-manifold with $%
\overset{0}{A}_{11}=0$. Then the solution of (\ref{2}) exists on $M\times
\lbrack 0,\infty ).$
\end{theorem}

There are many torsion free examples. On the other hand, this condition
implies strong topological obstruction (see the Appendix in \cite{ch}). We
hope this condition can be weakened in the future study.

To motivate the definitions of some analytic conditions, let us first
examine the standard $CR$ $3$-sphere $\mathbf{S}^{3}=\left\{
(z_{0},z_{1})|\sum_{j=0}^{1}z_{j}\bar{z}_{j}=1\right\} \subset \mathbf{C}%
^{2} $ with the induced $CR$ structure $J$ from $\mathbf{C}^{2}$ and the
contact form $\theta _{0}=\frac{i(\bar{\partial}u-\partial u)}{2}|_{\mathbf{S%
}^{3}}$, where $u=\left( \sum_{j=0}^{1}z_{j}\bar{z}_{j}\right) -1$ is a
defining function. With respect to $(J,\theta _{0})$, $\mathbf{S}^{3}$ is
torsion free and $Q_{0}=0$. On the other hand, if $f_{p,q}$ is a bigraded
spherical harmonic of type $(p,q)$ on $\mathbf{C}^{2}$ (i.e. a harmonic
polynomial which is a linear combination of terms of the form $z^{\rho }\bar{%
z}^{\gamma }$), then we have (\cite{chi}) ($\Delta _{0}$ denotes the
sub-Laplacian with respect to $(J,\theta _{0}))$%
\begin{equation}
\Delta _{0}f_{p,q}=-(2pq+p+q)\cdot f_{p,q}  \label{3}
\end{equation}%
and%
\begin{equation}
P_{0}f_{p,q}=4pq(p+1)(q+1)\cdot f_{p,q}.  \label{4}
\end{equation}

\begin{definition}
On a closed pseudohermitian manifold $(M,J,\theta ),$ we call the Paneitz
operator $P$ with respect to $(J,\theta )$ essentially positive if there
exists a constant $\Upsilon $ $>$ $0$ such that 
\begin{equation}
\int_{M}P\varphi \cdot \varphi d\mu \geq \Upsilon \int_{M}\varphi ^{2}d\mu .
\label{4'}
\end{equation}%
for all real $C^{\infty }$ smooth functions $\varphi $ $\bot \ Ker(P)$ (i.e.
perpendicular to the kernel of $P$ in the $L^{2}$ norm with respect to the
volume form $d\mu $ $=$ $\theta \wedge d\theta ).$
\end{definition}

\begin{remark}
\label{r0} We note that $P$ being essentially positive is a $CR$ invariant
property by (\ref{0a}), i.e., it is independent of the choice of contact
form. For $M$ being the boundary of a bounded strictly pseudoconvex domain
in $C^{2},$ $P$ appears in the transformation law of the first invariant in
the logarithmic term of Fefferman's asymptotic expansion of the Szeg\"{o}
kernel (\cite{h}). Also $P$ appears as the compatibility operator for the
degenerate Laplacian in the paper \cite{gl}. On the other hand, the kernel
of $P$ is infinite dimensional, containing all $CR$ -pluriharmonic functions
(see Section 5). So even for the short time solution to (\ref{2}), we need
to treat the kernel part separately. The condition that $P$ is essentially
positive was also used to study the problem of the first eigenvalue of the
sub-Laplacian (see \cite{chi}).
\end{remark}

Since the restrictions of bigraded spherical harmonics to $\mathbf{S}^{3}$
span a dense subspace of $L^{2}(\mathbf{S}^{3})$ (see Proposition 12.3.3 in 
\cite{cs}), we conclude that the Paneitz operator $P_{0}$ of $(\mathbf{S}%
^{3},J,\theta _{0})$ is essentially positive by (\ref{4}). Combining (\ref{3}%
) and (\ref{4}), we get 
\begin{equation*}
(2P_{0}-\Delta _{0}^{2})f_{p,q}=\lambda (p,q)f_{p,q},
\end{equation*}%
where $\lambda (p,q)=p^{2}(4q-1)+q^{2}(4p-1)+6pq+4p^{2}q^{2}$. It is clear
that if $P_{0}f_{p,q}\neq 0$, i.e., $pq\neq 0$, then $\lambda (p,q)>0$. This
means that the operator $2P_{0}-\Delta _{0}^{2}$ is positive on the
orthogonal complement of the Kernel of $P_{0}$. Therefore if $\varphi $ is a
real $C^{\infty }$ smooth function such that $\varphi \bot \ Ker(P_{0})$,
then 
\begin{equation}
2\int_{\mathbf{S}^{3}}P_{0}\varphi \cdot \varphi d\mu _{0}\geq \int_{\mathbf{%
S}^{3}}\Delta _{0}^{2}\varphi \cdot \varphi d\mu _{0}=\int_{\mathbf{S}%
^{3}}(\Delta _{0}\varphi )^{2}d\mu _{0}.  \label{5}
\end{equation}

Let $\overset{0}{\nabla }$ and $\overset{0}{\nabla ^{2}}$ denote the
sub-gradient and the sub-Hessian with respect to $(J,\theta _{0}),$
respectively. Now based on the Bochner formula, we have (Lemma \ref{l3.2} \
in Section $3;$ see \cite{chi} also) 
\begin{equation}
\begin{array}{ccl}
\int_{M}|\overset{0}{\nabla ^{2}}\varphi |^{2}d\mu _{0} & = & 
3\int_{M}(\Delta _{0}\varphi )^{2}d\mu _{0}-2\int_{M}P_{0}\varphi \cdot
\varphi d\mu _{0} \\ 
&  & -\int_{M}\overset{0}{W}|\overset{0}{\nabla }\varphi |^{2}d\mu _{0}-6%
\func{Im}\int_{M}\overset{0}{A}_{\overline{1}\overline{1}}\varphi
_{1}\varphi _{1}d\mu _{0}%
\end{array}
\label{6}
\end{equation}%
where $\overset{0}{W}$ and $\overset{0}{A}_{\overline{1}\overline{1}}$
denote the Tanaka-Webster curvature and the torsion with respect to $%
(J,\theta _{0}),$ respectively.

Observe that $\overset{0}{W}$ is a positive constant and $\overset{0}{A}_{%
\overline{1}\overline{1}}$ $=$ $0$ for $(\mathbf{S}^{3},J,\theta _{0})$. It
follows from (\ref{5}) and (\ref{6}) that 
\begin{equation}
\int_{\mathbf{S}^{3}}|\overset{0}{\nabla ^{2}}\varphi |^{2}d\mu _{0}\leq
2\int_{\mathbf{S}^{3}}(\Delta _{0}\varphi )^{2}d\mu _{0}  \label{7}
\end{equation}%
for $\varphi \bot \ Ker(P_{0}).$

Inspired by the inequality (\ref{7}), we make the following definition.

\begin{definition}
We say that the condition $(\ast )$ is satisfied on a closed CR $3$-manifold 
$(M,J,[\theta _{0}])$ if there exist constants $0\leq \varepsilon _{0}<1$
and $C(\varepsilon _{0})\geq 0$ such that \ 
\begin{equation}
\int_{M}|\overset{0}{\nabla ^{2}}\varphi |^{2}d\mu _{0}\leq (2+\varepsilon
_{0})\int_{M}(\Delta _{0}\varphi )^{2}d\mu _{0}+C(\varepsilon
_{0})\int_{M}\varphi ^{2}d\mu _{0}  \tag{$\ast $}
\end{equation}%
for all real $C^{\infty }$ smooth functions $\varphi $ $\bot \ Ker(P_{0}).$
\end{definition}

\begin{theorem}
\label{Th1.2} Let $(M,J,[\theta _{0}])$ be a closed CR $3$-manifold. Suppose
that $P_{0}$ is essentially positive and the condition $(\ast )$ holds. Then
the solution $\lambda $ of (\ref{2}) satisfies an a priori $S^{2,2}$ estimate%
$:$%
\begin{equation*}
||\lambda ||_{S^{2,2}}\leq C(T)
\end{equation*}%
for $t\in \lbrack 0,T),$ where $C(T)$ is a constant depending on $T.$
Moreover, if in addition $(Q_{0})_{\ker }=0,$ then the bound $C(T)$ can be
replaced by a constant $C$ independent of the time.
\end{theorem}

In the paper \cite{b}, the (essential) positivity of the Riemannian Paneitz
operator is also needed for a similar result on the $Q$-curvature flow on a
closed conformal $4$-manifold $(M,[g_{0}])$ with a given conformal class $%
[g_{0}]$. However, in the Riemannian case, the condition analogous to the
condition $(\ast )$ holds always. We wonder to what extent the condition $%
(\ast )$ is valid for a closed pseudohermitian manifold.

In Section $5$, we show that if the torsion $\overset{0}{A}_{11}$ of $%
(M,J,\theta _{0})$ is zero, then the $CR$ Paneitz operator $P_{0}$ is
essentially positive. In that case, the $CR$ Paneitz operator $P_{0}=\Box
_{b}\overline{\Box _{b}}$ and the Kohn Laplacian $\Box _{b}$ and $\overline{%
\Box _{b}}$ commute (see Section $5$). We can write an $L^{2}$ function $%
\varphi $ $=$ $\varphi _{\ker }$ $+$ $\varphi ^{\bot }$ where $\varphi
_{\ker }$ $\in $ $Ker(P_{0})$ and $\varphi ^{\bot }$ $\in $ $%
Ker(P_{0})^{\bot }.$ We have the following a priori estimates of higher
orders.

\begin{theorem}
\label{Th1.3} Let $(M,J,[\theta _{0}])$ be a closed CR $3$-manifold. Suppose
that $P_{0}$ is essentially positive and the condition $(\ast )$ holds. In
addition, suppose also that $\Delta _{0}(Ker(P_{0}))$ $\subset $ $%
Ker(P_{0}). $ Then for any nonnegative integer $k,$ the solution $\lambda $
of (\ref{2}) satisfies an a priori $S^{2k,2}$ estimate$:$%
\begin{equation*}
||\lambda ||_{S^{2k,2}}\leq C(k,T)
\end{equation*}%
for $t\in \lbrack 0,T),$ where $C(k,T)$ is a constant depending on $k$ and $%
T.$ Moreover, if in addition $(Q_{0})_{\ker }=0,$ then the bound $C(k,T)$
can be replaced by a constant $C(k)$ independent of the time.
\end{theorem}

We remark that in the torsion free case, the condition $\Delta
_{0}(Ker(P_{0}))$ $\subset $ $Ker(P_{0})$ holds true (also $P_{0}$ is
essentially positive as mentioned above). We have the following asymptotic
convergence of solutions of (\ref{2}).

\begin{theorem}
\label{t1} Let $(M,J,[\theta _{0}])$ be a closed CR $3$-manifold with $%
(Q_{0})_{\ker }=0$. Suppose that $\overset{0}{A}_{11}=0$. Then the solution
of (\ref{2}) exists on $M\times \lbrack 0,\infty )$ and converges smoothly
to $\lambda _{\infty }$ $\equiv $ $\lambda (\cdot ,\infty )$ as $t$ $%
\rightarrow $ $\infty .$ Moreover, the contact form $e^{2\lambda _{\infty
}}\theta _{0}$ has zero $Q$-curvature.
\end{theorem}

\begin{remark}
\label{r1} 1. Let \ $W_{0}$ denote the Tanaka-Webster (scalar) curvature
with respect to $(J,\theta _{0}).$ If $\overset{0}{A}_{11}=0$, then $P_{0}$
commutes with $\Delta _{0}$ and hence there holds 
\begin{equation*}
(Q_{0})_{\ker }=-\frac{2}{3}(\Delta _{0}W_{0})_{\ker }=-\frac{2}{3}\Delta
_{0}(W_{0})_{\ker }
\end{equation*}

by (\ref{-1}). It follows that 
\begin{equation*}
(Q_{0})_{\ker }=0\Longleftrightarrow (W_{0})_{\ker }\text{ is a constant.}
\end{equation*}

2. On the standard $CR$ $3$-sphere $(\mathbf{S}^{3},J,[\theta _{0}]),$ we
have

(i) $\overset{0}{A}_{11}=0$ and $Q_{0}=0$,

(ii) the condition $(\ast )$ holds with $\varepsilon _{0}=0$ and $%
C(\varepsilon _{0})=0$.
\end{remark}

As a consequence of Theorem \ref{t1}, we have

\begin{corollary}
\label{c1} Let $(\mathbf{S}^{3},J,[\theta _{0}])$ be the standard $CR$ $3$%
-sphere. Then the solution of (\ref{2}) exists on $\mathbf{S}^{3}\times
\lbrack 0,\infty )$ and converges smoothly to $\lambda _{\infty }$ such that 
$e^{2\lambda _{\infty }}\theta _{0}$ is a contact form of zero $Q$-curvature.
\end{corollary}

We recall that $\theta _{0}$ is called an invariant contact form on a $CR$ $3
$-manifold $M$ if it is locally volume- normalized with respect to a closed $%
(2,0)$-form on $M$ (\cite{fh}, \cite{l1}, \cite{fa}). In the paper \cite{fh}%
, the authors proved that the $Q$-curvature of an invariant contact form
vanishes. Indeed if $M$ is a real hypersurface in $\mathbf{C}^{2},$then $M$
admits an invariant contact form $\theta _{0}$ so that $Q_{0}=0$ on $M.$ In
general, there is a topological obstruction for the global existence of an
invariant contact form $\theta _{0}$ (\cite{l1}). However, on the s$\tan $%
dard $CR$ $3$-sphere$,$ $Q_{0}=0$ if and only if \ $\theta _{0}$ is an
invariant contact form.

About other curvature flows, we notice that it is still an open problem
whether we have the long-time existence and convergence for solutions of the 
$CR$ Yamabe flow on a closed $CR$ $3$-manifold (\cite{cc1}, \cite{jl}). This
is a second order subparabolic equation. On the other hand, the flow (\ref{2}%
) which we are dealing with is a fourth-order subparabolic (at least under a
certain condition) equation. Because there seems to be no suitable maximum
principle available for fourth-order subelliptic operators, we need to
invoke a priori $L^{2}$ estimates for solutions to (\ref{2}) in place of the
pointwise estimates used for second-order (sub)elliptic operators (\cite{b}, 
\cite{c}, \cite{c1}, \cite{c2}, \cite{cw}). In case $M$ is a surface, the $Q$%
-curvature flow corresponds to the $2$-dimensional Calabi flow which is
solved completely by P. T. Chru\'{s}ciel (\cite{c}) and the first author (%
\cite{c3}, \cite{c4}; see also \cite{cw}).

In order to get the $S^{2,2}$-estimate (see Appendix A for the definition of
Folland- Stein norms $S^{k,p})$, we need an additional analytic condition $%
(\ast )$ (which holds for the standard pseudohermitian $3$-sphere) plus a
trick from S. Brendle's \ work (\cite{b}) for the $Q$-curvature flow on
Riemannian $4$-manifolds. This is because the pseudohermitian version of a
Bochner-type estimate is entirely different from the Riemannian version (see
Lemma \ref{l3.1}).

We briefly describe the methods used in our proofs. In Section $2,$ in order
to have the $S^{k,2}$ -estimates for $\lambda ,$ we need to derive the
analogue of the Moser inequality on pseudohermitian $3$-manifolds. In
Section $3$, based on a pseudohermitian version of the Bochner formula, we
first derive a key estimate for the equation (\ref{2}) as in Lemma \ref{l3.2}%
, which involves the $CR$ Paneitz operator $P_{0}$. We show the
subellipticity of $P_{0}$ on ($\ker P_{0})^{\perp }$ under the torsion free
condition to get the short time solution. For the long-time solution of (\ref%
{2}) under the same condition (i.e., Theorem \ref{Th1.1}), we obtain the
higher order bounds for the solution.

In Section $4,$ we derive the $S^{2,2}$ -estimate and higher-order $S^{2k,2}$
-estimates for $\lambda $ under the flow (\ref{2}), and hence prove Theorem %
\ref{Th1.2} and Theorem \ref{Th1.3}. Then we prove the smooth convergence of
solutions of (\ref{2}) (i.e. Theorem \ref{t1}) by a method analogous to that
in \cite{s}. In Section $5,$ we show the essential positivity of the $CR$
Paneitz operator $P$ on $(M,J,\theta )$ with the zero torsion.

\section{\textbf{Moser's Inequality on Pseudohermitian }$3$\textbf{-manifolds%
}}

In this section, based on \cite{a}, \cite{cl}, \cite{fs1}, and \cite{sc}, we
derive the analogue of Moser's inequality on pseudohermitian $3$-manifolds.

Let $H^{1}=\mathbf{C}\times \mathbf{R}$\textbf{\ }be the ($3$-dimensional)
Heisenberg group with coordinates $(z,t)$. For each real number $r\in 
\mathbf{R}$, there is a dilation naturally associated with $H^{1},$ which is
usually denoted as 
\begin{equation*}
\delta _{r}u=\delta _{r}(z,t)=(rz,r^{2}t).
\end{equation*}

The anisotropic dilation structure on $H^{1}$ introduces a homogeneous norm%
\begin{equation*}
\left\vert u\right\vert =\left\vert (z,t)\right\vert =\left( \left\vert
z\right\vert ^{4}+t^{2}\right) ^{1/4}.
\end{equation*}
With this norm, we can define the Heisenberg ball centered at $u=(z,t)$ with
radius $R$ by $B(u,R)=\left\{ v\in H^{1}:\left\vert u^{-1}\cdot v\right\vert
<R\right\} .$

In their paper (\cite{cl}), William S. Cohn and Guozhen Lu show that for all 
$\varphi \in C_{0}^{\infty }\left( H^{1}\right) $,

\begin{equation}
\left| \varphi (v)\right| \leq L^{-1}\int_{H^{1}}\frac{\left| \nabla
_{b}\varphi (u)\right| }{\left| v^{-1}\cdot u\right| ^{3}}\ dV(u),
\label{ineq1}
\end{equation}%
where $L=\frac{2\pi \ \Gamma (1/2)\Gamma (3/4)}{\Gamma (1)\Gamma (5/4)}$ and 
$dV(u)=dx\wedge dy\wedge dt,$ $z=x+iy$. Note that we only take the case $%
\beta =1$ in Theorem 1.2 of (\cite{cl}).

Let $B$ be the unit Heisenberg ball of $H^{1},$ centered at $(0,0)$. Let $%
||\cdot ||_{p}$ denote the $L^{p}$ norm with respect to the volume form $%
dV(u).$

\begin{lemma}
For $\varphi \in C_{0}^{\infty }(B)$ and $p\geq 4$, we have 
\begin{equation}
\left\| \varphi \right\| _{p}\leq L^{-1}\left\| \nabla _{b}\varphi \right\|
_{4}\ \sup_{v\in B}\left[ \int_{B}\left| v^{-1}\cdot u\right| ^{-3k}\ dV(u)%
\right] ^{\frac{1}{k}}  \label{inequ2}
\end{equation}
where $\frac{1}{k}=\frac{1}{p}+\frac{3}{4}.$
\end{lemma}

\proof
\ We write

\begin{equation*}
\begin{split}
& \left\vert v^{-1}\cdot u\right\vert ^{-3}|\nabla _{b}\varphi (u)| \\
=& \left( \left\vert v^{-1}\cdot u\right\vert ^{-3k}\left\vert \nabla
_{b}\varphi (u)\right\vert ^{4}\right) ^{\frac{1}{p}}\left( \left\vert
v^{-1}\cdot u\right\vert ^{-3k}\right) ^{\frac{3}{4}}\left( |\nabla
_{b}\varphi (u)|^{4}\right) ^{\frac{1}{4}-\frac{1}{p}}.
\end{split}%
\end{equation*}%
Since $\frac{1}{p}+\frac{1}{4}+\left( \frac{1}{4}-\frac{1}{p}\right) =1$,
applying H\"{o}lder's inequality to (\ref{ineq1}), we have

\begin{equation*}
\begin{split}
\left\vert \varphi (v)\right\vert \leq L^{-1}& \left( \int_{B}\left\vert
v^{-1}\cdot u\right\vert ^{-3k}|\nabla _{b}\varphi (u)|^{4}\ dV(u)\right) ^{%
\frac{1}{p}}\left( \int_{B}\left\vert v^{-1}\cdot u\right\vert ^{-3k}\
dV(u)\right) ^{\frac{3}{4}} \\
& \left( \int_{B}|\nabla _{b}\varphi (u)|^{4}\ dV(u)\right) ^{\frac{1}{4}-%
\frac{1}{p}}.
\end{split}%
\end{equation*}%
This implies that 
\begin{equation*}
\begin{split}
\left\Vert \varphi \right\Vert _{p}& =\left( \int_{B}\left\vert \varphi
(v)\right\vert ^{p}\ dV(v)\right) ^{\frac{1}{p}} \\
& \leq \left( L^{-1}\left\Vert \nabla _{b}\varphi \right\Vert _{4}^{1-\frac{4%
}{p}}\cdot \Lambda ^{\frac{3k}{4}}\right) \left[ \int_{B}\int_{B}\left\vert
v^{-1}\cdot u\right\vert ^{-3k}|\nabla _{b}\varphi (u)|^{4}\ dV(u)dV(v)%
\right] ^{\frac{1}{p}} \\
& \leq L^{-1}\left\Vert \nabla _{b}\varphi \right\Vert _{4}\cdot \Lambda ,
\end{split}%
\end{equation*}%
where $\Lambda =\sup_{v\in B}\left[ \int_{B}\left\vert v^{-1}\cdot
u\right\vert ^{-3k}\ dV(u)\right] ^{\frac{1}{k}}$.

\endproof%

\begin{corollary}
For all $\varphi \in C_{0}^{\infty }(B)$, 
\begin{equation}
\int_{B}e^{\varphi }\ dV(u)\leq C\exp {\left( \varkappa \left\| \nabla
_{b}\varphi \right\| _{4}^{4}\right) },  \label{ineq3}
\end{equation}%
where $C$ and $\varkappa $ are two positive constants.
\end{corollary}

\proof
Let $\Sigma =\partial B$ be the unit Heisenberg sphere and $dA$ be the
unique Radon measure on $\Sigma $ (\cite{cl}, \cite{fs1}). We denote $%
A_{1}=\int_{\Sigma }\ dA$. Note that there exists a number $\delta $ such
that $\left\vert v^{-1}\cdot u\right\vert <\delta $ for all $u,v\in B$. This
means that $B\subset B(v,\delta )$ for all $v\in B$.

By (\ref{inequ2}), we have 
\begin{equation}
\begin{split}
\left\Vert \varphi \right\Vert _{p}& \leq L^{-1}\left\Vert \nabla
_{b}\varphi \right\Vert _{4}\ \sup_{v\in B}\left[ \int_{B}\left\vert
v^{-1}\cdot u\right\vert ^{-3k}\ dV(u)\right] ^{\frac{1}{k}} \\
& \leq L^{-1}\left\Vert \nabla _{b}\varphi \right\Vert _{4}\ \sup_{v\in B} 
\left[ \int_{B(v,\delta )}\left\vert v^{-1}\cdot u\right\vert ^{-3k}\
dV(v^{-1}\cdot u)\right] ^{\frac{1}{k}} \\
& =L^{-1}\left\Vert \nabla _{b}\varphi \right\Vert _{4}\ \left[
\int_{B(0,\delta )}\left\vert u\right\vert ^{-3k}\ dV(u)\right] ^{\frac{1}{k}%
} \\
& =L^{-1}A_{1}^{\frac{1}{k}}\left\Vert \nabla _{b}\varphi \right\Vert
_{4}\left( \int_{0}^{\delta }\ r^{3-3k}\ dr\right) ^{\frac{1}{k}} \\
& =L^{-1}A_{1}^{\frac{1}{k}}\delta ^{\frac{4-3k}{k}}\left\Vert \nabla
_{b}\varphi \right\Vert _{4}\left( \frac{1}{4-3k}\right) ^{\frac{1}{k}}.
\end{split}
\label{ineq4}
\end{equation}%
Since $\frac{1}{k}=\frac{1}{p}+\frac{3}{4}$ and $p\geq 4$, we immediately
get that $\left( \frac{1}{4-3k}\right) ^{\frac{1}{k}}\leq Cp^{\frac{3}{4}}$
and $-3\leq \frac{4-3k}{k}\leq 1$ for some constant $C$. Thus, by (\ref%
{ineq4}), there exists a constant $K$ such that for all $p\geq 1$,

\begin{equation}
\left\| \varphi \right\| _{p}\leq K\left\| \nabla _{b}\varphi \right\|
_{4}\cdot p^{\frac{3}{4}}.  \label{ineq5}
\end{equation}%
It follows that

\begin{equation}
\begin{split}
\int_{B}\ e^{\varphi }\ dV(u)& =\int_{B}\ \left( \sum_{p=0}^{\infty }\frac{%
\varphi ^{p}}{p!}\right) \ dV(u) \\
& \leq \sum_{p=0}^{\infty }\frac{\left\Vert \varphi \right\Vert _{p}^{p}}{p!}
\\
& \leq \sum_{p=0}^{\infty }K^{p}\left\Vert \nabla _{b}\varphi \right\Vert
_{4}^{p}(p!)^{-1}p^{\frac{3p}{4}} \\
& =\sum_{p=0}^{\infty }\frac{\left( K^{4}\right) ^{\frac{p}{4}}\left(
\left\Vert \nabla _{b}\varphi \right\Vert _{4}^{4}\right) ^{\frac{p}{4}}}{(%
\frac{p}{4})!}\frac{(\frac{p}{4})!}{p!}\ p^{\frac{3p}{4}}.
\end{split}
\label{ineq6}
\end{equation}%
Here $x!=\Gamma (x+1)$ for all real number $x\geq 0$.

According to Stirling's formula, when $p\rightarrow \infty $, we can estimate

\begin{equation*}
\begin{split}
\frac{(\frac{p}{4})!}{p!}\ p^{\frac{3p}{4}}& \approx \frac{\left( 2\pi \frac{%
p}{4}\right) ^{\frac{1}{2}}\left( \frac{p}{4e}\right) ^{\frac{p}{4}}}{\left(
2\pi p\right) ^{\frac{1}{2}}\left( \frac{p}{e}\right) ^{p}}\ p^{\frac{3p}{4}}
\\
& \approx \frac{1}{2}\left( \frac{1}{4}\right) ^{\frac{p}{4}}e^{\frac{3p}{4}%
}.
\end{split}%
\end{equation*}%
Therefore, from (\ref{ineq6}), we get

\begin{equation*}
\begin{split}
\int_{B}\ e^{\varphi }\ dV(u)& \leq C\sum_{p=0}^{\infty }\frac{\left( \frac{%
e^{3}}{4}K^{4}\right) ^{\frac{p}{4}}\left( \left\| \nabla _{b}\varphi
\right\| _{4}^{4}\right) ^{\frac{p}{4}}}{(\frac{p}{4})!} \\
& \leq C\cdot \exp {\left( \varkappa \left\| \nabla _{b}\varphi \right\|
_{4}^{4}\right) },
\end{split}%
\end{equation*}%
for some constants $C$ and $\varkappa $.

\endproof%

Now we are ready to prove an analogue of the Moser inequality in
pseudohermitian geometry. Let $C^{\infty }(M)$ denote the space of all real
valued $C^{\infty }$ smooth functions on $M.$

\begin{theorem}
\label{T2.1} (Pseudohermitian Moser inequality) Let $(M,J,\theta )$ be a
closed pseudohermitian $3$-manifold. Then there exist constants $C,$ $%
\varkappa ,$ and $\nu $ such that for all $\varphi \in C^{\infty }(M)$,
there holds 
\begin{equation}
\int_{M}e^{\varphi }\ d\mu \leq C\exp {\left( \varkappa \left\Vert \nabla
_{b}\varphi \right\Vert _{4}^{4}+\nu \left\Vert \varphi \right\Vert
_{4}^{4}\right) }  \label{moine}
\end{equation}%
where $d\mu =\theta \wedge d\theta $ and the $L^{4}$ norm $||\cdot ||_{4}$
respects the volume form $d\mu $.
\end{theorem}

\proof
For each point $x\in M$, there exists a neighborhood $U_{x}$ of $x$ such
that $U_{x}$ is diffemorphic to a Heisenberg ball $B(r_{x})$ centered at the
origin with radius $r_{x}<1$. Choose a cut-off function $\eta _{x}$ such
that $\eta _{x}(v)=1$ for $\left| v\right| <\frac{r_{x}}{2}$ and $\eta
_{x}(v)=0$ for $\left| v\right| \geq r_{x}$.

Since $M$ is closed, there exist finite balls $\left( B(\frac{r_{j}}{2}%
),\eta _{j}\right) ,\ j=1,\cdots m$ such that $\left\{ B(\frac{r_{j}}{2}%
)\right\} $ is an open covering of $M$. Let $\varphi _{j}=\eta _{j}\varphi $%
. Making use of (\ref{ineq3}), we compute 
\begin{equation*}
\begin{split}
\int_{M}\ e^{\varphi }\ d\mu & \leq C\cdot \sum_{j=1}^{m}\int_{B(\frac{r_{j}%
}{2})}\ e^{\varphi }\ dV(u) \\
& \leq C\cdot \sum_{j=1}^{m}\int_{B}\ e^{\varphi _{j}}\ dV(u) \\
& \leq C\cdot \sum_{j=1}^{m}\exp {\left( \varkappa \left\| \nabla
_{b}\varphi _{j}\right\| _{4;B}^{4}\right) } \\
& \leq C\cdot \sum_{j=1}^{m}\exp {\left( \varkappa \left\| \nabla
_{b}\varphi \right\| _{4}^{4}+\nu \left\| \varphi \right\| _{4}^{4}\right) },
\end{split}%
\end{equation*}%
for some constants $C,\ \varkappa ,$ and $\nu $.

\endproof%

\section{\textbf{The Long-time Existence}}

Let $T$ be the maximal time for a solution of the flow (\ref{2}) on $M\times
\lbrack 0,T)$. We will derive the $S^{k,2}$-norm estimate for $\lambda $
under the flow (\ref{2}) for all $0\leq t<T$. It then follows that we have
the long-time existence for solutions of (\ref{2}) on $M\times \lbrack
0,\infty ).$

First we have an integral version of the\ Bochner formula on a
pseudohermitian $3$-manifold.

\begin{lemma}
\label{l3.1}Let $(M,J,\theta _{0})$ be a closed pseudohermitian $3$%
-manifold. Then for any $\lambda $ $\in $ $C^{\infty }(M),$ there holds%
\begin{equation}
\begin{array}{cl}
0= & \int_{M}(\Delta _{0}\lambda )^{2}d\mu _{0}-\int_{M}|\overset{0}{\nabla
^{2}}\lambda |^{2}d\mu _{0}+2\int_{M}(\lambda _{0})^{2}d\mu _{0} \\ 
& -\int_{M}\overset{0}{W}|\overset{0}{\nabla }\lambda |^{2}d\mu _{0}+2\func{%
Im}\int_{M}\overset{0}{A}_{\overline{1}\overline{1}}\lambda _{1}\lambda
_{1}d\mu _{0}%
\end{array}
\label{3.1}
\end{equation}%
where $\overset{0}{\nabla }$ and $T_{0}$ denote the sub-gradient and the
characteristic vector field with respect to $(J,\theta _{0}),$ respectively,
and $\lambda _{0}=T_{0}\lambda .$
\end{lemma}

\proof
We first show that (see Appendix A for definitions of the notations)%
\begin{equation}
\begin{array}{ccl}
\frac{1}{2}\Delta _{0}|\overset{0}{\nabla }\lambda |^{2} & = & |\overset{0}{%
\nabla ^{2}}\lambda |^{2}+<\overset{0}{\nabla }\lambda ,\overset{0}{\nabla }%
(\Delta _{0}\lambda )>_{J,\theta _{0}}+\overset{0}{W}|\overset{0}{\nabla }%
\lambda |^{2} \\ 
&  & +Tor(\overset{0}{\nabla }\lambda ,\overset{0}{\nabla }\lambda
)-2i\lambda _{1}\lambda _{0\overline{1}}+2i\lambda _{\overline{1}}\lambda
_{01},%
\end{array}
\label{10}
\end{equation}%
where $Tor(\overset{0}{\nabla }\lambda ,\overset{0}{\nabla }\lambda )\equiv i%
\overset{0}{A}_{\overline{1}\overline{1}}\lambda _{1}\lambda _{1}+\mathrm{%
conjugate\ of}$ $i\overset{0}{A}_{\overline{1}\overline{1}}\lambda
_{1}\lambda _{1}.$

We compute 
\begin{equation*}
|\overset{0}{\nabla }\lambda |^{2}=2|\lambda _{1}|^{2}=2\lambda _{1}\lambda
_{\overline{1}},\ |\overset{0}{\nabla ^{2}}\lambda |^{2}=2\lambda _{%
\overline{1}\overline{1}}\lambda _{11}+2\lambda _{\overline{1}1}\lambda _{1%
\overline{1}}
\end{equation*}%
and 
\begin{equation}
\begin{array}{ccl}
\frac{1}{2}\Delta _{0}|\overset{0}{\nabla }\lambda |^{2} & = & (\lambda
_{1}\lambda _{\overline{1}})_{1\overline{1}}+(\lambda _{1}\lambda _{%
\overline{1}})_{\overline{1}1} \\ 
& = & \lambda _{\overline{1}1\overline{1}}\lambda _{1}+2\lambda _{\overline{1%
}1}\lambda _{1\overline{1}}+2\lambda _{\overline{1}\overline{1}}\lambda
_{11}+\lambda _{\overline{1}}\lambda _{11\overline{1}}+\lambda _{1}\lambda _{%
\overline{1}\overline{1}1}+\lambda _{\overline{1}}\lambda _{1\overline{1}1}
\\ 
& = & |\overset{0}{\nabla ^{2}}\lambda |^{2}+\lambda _{1}(\lambda _{%
\overline{1}\overline{1}1}+\lambda _{\overline{1}1\overline{1}})+\mathrm{%
conjugate\ of}\text{\ }\lambda _{1}(\lambda _{\overline{1}\overline{1}%
1}+\lambda _{\overline{1}1\overline{1}}).%
\end{array}
\label{12}
\end{equation}%
Observe that 
\begin{equation}
\begin{array}{ccc}
<\overset{0}{\nabla }\lambda ,\overset{0}{\nabla }(\Delta _{0}\lambda
)>_{J,\theta _{0}} & = & \lambda _{1}(\lambda _{1\overline{1}}+\lambda _{%
\overline{1}1})_{\overline{1}}+\mathrm{conjugate\ of\ }\lambda _{1}(\lambda
_{1\overline{1}}+\lambda _{\overline{1}1})_{\overline{1}} \\ 
& = & \lambda _{1}(\lambda _{1\overline{1}\overline{1}}+\lambda _{\overline{1%
}1\overline{1}})+\mathrm{conjugate\ of\ }\lambda _{1}(\lambda _{1\overline{1}%
\overline{1}}+\lambda _{\overline{1}1\overline{1}}).%
\end{array}
\label{14}
\end{equation}
\ 

Taking the covariant differentiation of%
\begin{equation}
\lambda _{\overline{1}1}=\lambda _{1\overline{1}}-i\lambda _{0}  \label{19}
\end{equation}%
(\cite{l1}) in the $Z_{\bar{1}}$ direction gives 
\begin{equation}
\begin{array}{ccc}
\lambda _{\overline{1}1\overline{1}} & = & \lambda _{1\overline{1}\overline{1%
}}-i\lambda _{0\overline{1}}.%
\end{array}
\label{16}
\end{equation}%
A commutation relation (\cite{l1}) for covariant derivatives of a 1-form
gives%
\begin{equation}
\lambda _{\overline{1}\overline{1}1}=\lambda _{\overline{1}1\overline{1}%
}-i\lambda _{\overline{1}0}+\overset{0}{W}\lambda _{\overline{1}}.
\label{18}
\end{equation}%
Now from (\ref{12}), (\ref{14}), (\ref{16}), and (\ref{18}), we obtain%
\begin{equation}
\begin{array}{ccl}
\frac{1}{2}\Delta _{0}|\overset{0}{\nabla }\lambda |^{2} & = & |\overset{0}{%
\nabla ^{2}}\lambda |^{2}+<\overset{0}{\nabla }\lambda ,\overset{0}{\nabla }%
(\Delta _{0}\lambda )>+(\overset{0}{W}\lambda _{1}\lambda _{\overline{1}}+%
\mathrm{conjugate}) \\ 
&  & +[-i\lambda _{1}(\lambda _{\overline{1}0}+\lambda _{0\overline{1}})+%
\mathrm{conjugate}].%
\end{array}
\label{20}
\end{equation}%
Then (\ref{10}) follows from (\ref{20}) and the following commutation
relation (\cite{l1})%
\begin{equation*}
\lambda _{0\overline{1}}=\lambda _{\overline{1}0}+\overset{0}{A}_{\overline{1%
}\overline{1}}\lambda _{1}.
\end{equation*}

Finally integrating both sides of (\ref{10}) and applying (\ref{19}), we
obtain%
\begin{equation*}
\begin{array}{cl}
0= & \int_{M}(\Delta _{0}\lambda )^{2}d\mu _{0}-\int_{M}|\overset{0}{\nabla
^{2}}\lambda |^{2}d\mu _{0}+2\int_{M}(i\lambda _{1}\lambda _{0\overline{1}%
}-i\lambda _{\overline{1}}\lambda _{01})d\mu _{0} \\ 
& -\int_{M}\overset{0}{W}|\overset{0}{\nabla }\lambda |^{2}d\mu _{0}+2\func{%
Im}\int_{M}\overset{0}{A}_{\overline{1}\overline{1}}\lambda _{1}\lambda
_{1}d\mu _{0}.%
\end{array}%
\end{equation*}%
By integrating by parts and (\ref{19}), we compute%
\begin{eqnarray*}
\int_{M}(i\lambda _{1}\lambda _{0\overline{1}}-i\lambda _{\overline{1}%
}\lambda _{01})d\mu _{0} &=&\int_{M}i\lambda _{0}(-\lambda _{1\overline{1}%
}+\lambda _{\overline{1}1})d\mu _{0} \\
&=&\int_{M}(\lambda _{0})^{2}d\mu _{0}.
\end{eqnarray*}%
Combining the above two formulas gives (\ref{3.1}).

\endproof%

\begin{lemma}
\label{l3.2} Let $(M,J,\theta _{0})$ be a closed pseudohermitian $3$%
-manifold. Then for any $\lambda $ $\in $ $C^{\infty }(M),$ there holds%
\begin{equation*}
\begin{array}{ccl}
2\int_{M}P_{0}\lambda \cdot \lambda d\mu _{0} & = & 3\int_{M}(\Delta
_{0}\lambda )^{2}d\mu _{0}-\int_{M}|\overset{0}{\nabla ^{2}}\lambda
|^{2}d\mu _{0} \\ 
&  & -\int_{M}\overset{0}{W}|\overset{0}{\nabla }\lambda |^{2}d\mu _{0}-6%
\func{Im}\int_{M}\overset{0}{A}_{\overline{1}\overline{1}}\lambda
_{1}\lambda _{1}d\mu _{0}.%
\end{array}%
\end{equation*}
\end{lemma}

\proof
Multiplying both sides of the formula $P_{0}\lambda =\Delta _{0}^{2}\lambda
+T_{0}^{2}\lambda +4\func{Im}(\overset{0}{A}_{\overline{11}}\lambda _{11}+%
\overset{0}{A}_{\overline{11},1}\lambda _{1})$ by $\lambda $ and
integrating, we compute

\begin{equation*}
\begin{array}{ccl}
2\int_{M}P_{0}\lambda \cdot \lambda d\mu _{0} & = & 2\int_{M}(\Delta
_{0}\lambda )^{2}d\mu _{0}-2\int_{M}(\lambda _{0})^{2}d\mu _{0}-8\func{Im}%
\int_{M}\overset{0}{A}_{\overline{1}\overline{1}}\lambda _{1}\lambda
_{1}d\mu _{0} \\ 
& = & 3\int_{M}(\Delta _{0}\lambda )^{2}d\mu _{0}-\int_{M}|\overset{0}{%
\nabla ^{2}}\lambda |^{2}d\mu _{0}-\int_{M}\overset{0}{W}|\overset{0}{\nabla 
}\lambda |^{2}d\mu _{0} \\ 
&  & -6\func{Im}\int_{M}\overset{0}{A}_{\overline{1}\overline{1}}\lambda
_{1}\lambda _{1}d\mu _{0}%
\end{array}%
\end{equation*}

by integrating by parts and Lemma \ref{l3.1}.

\endproof%

As a consequence of Lemma \ref{l3.2} and (\ref{5}), we have

\begin{corollary}
\label{c3.1}Let $(\mathbf{S}^{3},J,\theta _{0})$ be the standard
pseudohermitian $3$-sphere$.$ Then 
\begin{equation*}
\int_{\mathbf{S}^{3}}|\overset{0}{\nabla ^{2}}\lambda |^{2}d\mu _{0}\leq
2\int_{\mathbf{S}^{3}}(\Delta _{0}\lambda )^{2}d\mu _{0}
\end{equation*}%
for $\lambda \perp $ker$(P_{0}).$
\end{corollary}

It follows from (\ref{2d}) that

\begin{lemma}
\label{l3.3} Let $(M,J,[\theta _{0}])$ be a closed CR $3$-manifold. Let $%
\lambda $ be a solution of the flow (\ref{2}) on $M\times \lbrack 0,T)$.
Then there exists a positive constant $\beta =\beta (Q_{0},\theta _{0})$
such that%
\begin{equation}
\mathcal{E}\emph{(}\theta \emph{)}=\int_{M}P_{0}\lambda \cdot \lambda d\mu
_{0}+\int_{M}Q_{0}\lambda d\mu _{0}\leq \beta ^{2},  \label{9}
\end{equation}%
for all t $\in $ $[0,T).$
\end{lemma}

\begin{lemma}
\label{ODE} Let $f$ $:[0,T)\rightarrow R$ be a $C^{1}$ smooth function
satisfying $f^{\prime }$ $\leq $ $-C_{1}f+C_{2}$ for some positive constants 
$C_{1},$ $C_{2}>0.$ Then $f(t)$ $\leq $ $f(0)e^{-C_{1}t}+\frac{C_{2}}{C_{1}}$
for $t$ $\in $ $[0,T).$
\end{lemma}

We will often use the above lemma (whose proof is left to the reader) to
obtain the higher order estimates. Now we write $\lambda =\lambda _{\ker
}+\lambda ^{\perp }.$ $Q_{0}=(Q_{0})_{\ker }+Q_{0}^{\perp }$ with respect to 
$P_{0}.$ Comparing both sides of the following formula:

\begin{equation*}
\frac{\partial \lambda _{\ker }}{\partial t}+\frac{\partial \lambda ^{\perp }%
}{\partial t}=\frac{\partial \lambda }{\partial t}=-(Q_{0}+2P_{0}\lambda
)+r(t).
\end{equation*}%
we obtain%
\begin{equation}
\frac{\partial \lambda ^{\perp }}{\partial t}=-(Q_{0}^{\perp }+2P_{0}\lambda
^{\perp })  \label{11}
\end{equation}%
and 
\begin{equation}
\frac{\partial \lambda _{\ker }}{\partial t}=-(Q_{0})_{\ker }+r(t).
\label{11a}
\end{equation}

From now on, $C$ or $C_{j\text{ }}$($C(k,T),$ $C(T),$ etc., respectively)
denotes a generic constant (with emphasis on depending on $(k,T),$ $T,$
etc., respectively) which may vary from line to line.

\begin{proposition}
\label{p3.1} Let $(M,J,[\theta _{0}])$ be a closed CR $3$-manifold with $%
\overset{0}{A}_{11}=0$. \ Then under the flow (\ref{2}) (or equivalently (%
\ref{11}) and (\ref{11a})), there exists a positive constant $C(k,\Upsilon
,\theta _{0},\beta ,T)>0$ such that%
\begin{equation*}
||\lambda ||_{S^{2k,2}}\leq C(k,\Upsilon ,\theta _{0},\beta ,T)
\end{equation*}%
for all $0\leq t<T.$
\end{proposition}

\proof
From (\ref{11a}), we have 
\begin{equation}
\lambda _{\ker }(t,p)-\lambda _{\ker }(0,p)=-(Q_{0})_{\ker
}(p)t+\int_{0}^{t}rdt.  \label{1.1}
\end{equation}%
Since ($P_{0}$ being self-adjoint)%
\begin{equation*}
\int_{M}P_{0}\lambda ^{\perp }\cdot \lambda _{\ker }d\mu _{0}=0\text{ }
\end{equation*}%
and 
\begin{equation*}
\begin{array}{l}
\int_{M}Q_{0}\lambda _{\ker }d\mu _{0} \\ 
=\int_{M}Q_{0}[\lambda _{\ker }(0,p)-(Q_{0})_{\ker
}(p)t+\int_{0}^{t}rdt]d\mu _{0} \\ 
=\int_{M}Q_{0}\lambda _{\ker }(0,p)d\mu _{0}-t\int_{M}[(Q_{0})_{\ker
}(p)]^{2}d\mu _{0} \\ 
\geq -C-Ct,%
\end{array}%
\end{equation*}%
we compute%
\begin{equation*}
\begin{array}{l}
\int_{M}P_{0}\lambda \cdot \lambda d\mu _{0}+\int_{M}Q_{0}\lambda d\mu _{0}
\\ 
=\int_{M}P_{0}\lambda ^{\perp }\cdot (\lambda _{\ker }+\lambda ^{\perp
})d\mu _{0}+\int_{M}Q_{0}(\lambda _{\ker }+\lambda ^{\perp })d\mu _{0} \\ 
\geq \int_{M}P_{0}\lambda ^{\perp }\cdot \lambda ^{\perp }d\mu
_{0}+\int_{M}Q_{0}\lambda ^{\perp }d\mu _{0}-C-Ct.%
\end{array}%
\end{equation*}%
It then follows from Lemma \ref{l3.3} that 
\begin{equation*}
\int_{M}P_{0}\lambda ^{\perp }\cdot \lambda ^{\perp }d\mu
_{0}+\int_{M}Q_{0}\lambda ^{\perp }d\mu _{0}\leq (\beta ^{2}+C)+Ct
\end{equation*}%
for all $t$ $\in $ $[0,T).$ Since the torsion of $(M,J,\theta _{0})$
vanishes, the CR Paneitz operator $P_{0}$ is essentially positive in the
sense that 
\begin{equation*}
\int_{M}P_{0}\lambda ^{\perp }\cdot \lambda ^{\perp }d\mu _{0}\geq \Upsilon
\int_{M}(\lambda ^{\perp })^{2}d\mu _{0}
\end{equation*}%
by Theorem \ref{t2.1}. All these imply that there exists a positive constant 
$C(\Upsilon ,Q_{0},\theta _{0})>0$ such that%
\begin{equation*}
\begin{array}{ccc}
(\beta ^{2}+C)+Ct & \geq & \int_{M}P_{0}\lambda ^{\perp }\cdot \lambda
^{\perp }d\mu _{0}+\int_{M}Q_{0}\lambda ^{\perp }d\mu _{0} \\ 
& \geq & \frac{\Upsilon }{2}\int_{M}(\lambda ^{\perp })^{2}d\mu
_{0}-C(\Upsilon ,Q_{0},\theta _{0})%
\end{array}%
\end{equation*}%
for all $t$ $\in $ $[0,T).$ So there exists a positive constant $C(\Upsilon
,Q_{0},\theta _{0},\beta ,T)$ such that 
\begin{equation}
\int_{M}(\lambda ^{\perp })^{2}d\mu _{0}\leq C(\Upsilon ,Q_{0},\theta
_{0},\beta ,T),  \label{1.2}
\end{equation}%
for all $t$ $\in $ $[0,T).$

On the other hand, we observe from the condition $\overset{0}{A}_{11}=0$ that

\begin{equation}
\Delta _{0}^{k}T_{0}^{2}=T_{0}^{2}\Delta _{0}^{k}  \label{1.5}
\end{equation}%
and hence

\begin{equation}
\Delta _{0}^{k}P_{0}=P_{0}\Delta _{0}^{k}.  \label{1.5'}
\end{equation}
It follows that

\begin{equation}
\Delta _{0}^{k}\lambda ^{\perp }\perp \ker P_{0},\ \Delta _{0}^{k}\lambda
_{\ker }\in \ker P_{0}  \label{1.4}
\end{equation}

Next we compute, for all positive integers $k,$ 
\begin{equation}
\begin{array}{l}
\frac{d}{dt}\int_{M}(\Delta _{0}^{k}\lambda ^{\perp })^{2}d\mu _{0} \\ 
=2\int_{M}(\Delta _{0}^{k}\lambda ^{\perp })(\Delta _{0}^{k}\frac{\partial
\lambda ^{\perp }}{\partial t})d\mu _{0} \\ 
=-2\int_{M}(\Delta _{0}^{k}\lambda ^{\perp })[\Delta _{0}^{k}(Q_{0}^{\perp
}+2P_{0}\lambda ^{\perp })]d\mu _{0} \\ 
=-2\int_{M}\lambda ^{\perp }[\Delta _{0}^{2k}(Q_{0}^{\perp })]d\mu
_{0}-4\int_{M}(\Delta _{0}^{k}\lambda ^{\perp })(P_{0}\Delta _{0}^{k}\lambda
^{\perp })d\mu _{0}.%
\end{array}
\label{eqn1.6}
\end{equation}%
Here we have used (\ref{1.5'}). By (\ref{1.4}) and essential positivity of $%
P_{0},$ we obtain%
\begin{equation}
\int_{M}(\Delta _{0}^{k}\lambda ^{\perp })(P_{0}\Delta _{0}^{k}\lambda
^{\perp })d\mu _{0}\geq \Upsilon \int_{M}(\Delta _{0}^{k}\lambda ^{\perp
})^{2}d\mu _{0}.  \label{eqn1.7}
\end{equation}%
Therefore from (\ref{eqn1.6}), (\ref{eqn1.7}), and (\ref{1.2}), we conclude
that%
\begin{equation}
\frac{d}{dt}\int_{M}(\Delta _{0}^{k}\lambda ^{\perp })^{2}d\mu _{0}\leq
-4\Upsilon \int_{M}(\Delta _{0}^{k}\lambda ^{\perp })^{2}d\mu _{0}+C(k,T).
\label{eqn1.7'}
\end{equation}%
By applying Lemma \ref{ODE} to the O.D.E. $f^{\prime }(t)\leq -4\Upsilon
f(t)+C(k,T)$ from (\ref{eqn1.7'}), we obtain 
\begin{equation*}
\int_{M}(\Delta _{0}^{k}\lambda ^{\perp })^{2}d\mu _{0}\leq C(k,T).
\end{equation*}%
and hence%
\begin{equation*}
||\lambda ^{\perp }-\overline{\lambda ^{\perp }}||_{S^{2k,2}}\leq C(k,T)
\end{equation*}%
for all $0\leq t<T.$

From the definition of the $Q$-curvature, $\overset{0}{A}_{11}=0,$ and (\ref%
{1.4}), we compute 
\begin{equation*}
\int_{M}(Q_{0})_{\ker }d\mu _{0}=-\frac{2}{3}\int_{M}\Delta
_{0}(R_{0})_{\ker }d\mu _{0}=0.
\end{equation*}%
This and (\ref{1.1}) imply%
\begin{equation*}
\lambda -\overline{\lambda }=(\lambda ^{\perp }-\overline{\lambda ^{\perp }}%
)+[\lambda _{\ker }(0,p)-(Q_{0})_{\ker }(p)t]-\overline{\lambda _{\ker }(0,p)%
},
\end{equation*}%
and hence 
\begin{equation}
||\lambda -\overline{\lambda }||_{S^{2k,2}}\leq C((Q_{0})_{\ker
},k,T)+||\lambda ^{\perp }-\overline{\lambda ^{\perp }}||_{S^{2k,2}}\leq
C(k,T)  \label{eqn1.8}
\end{equation}%
for all $0\leq t<T.$ Recall that the average $\bar{f}$ of a function $f$ is
defined by $\overline{f}=\frac{\int_{M}fd\mu _{0}}{\int_{M}d\mu _{0}}.$ In
particular, there holds%
\begin{equation*}
||\lambda -\overline{\lambda }||_{S^{2,2}}\leq C(T).
\end{equation*}%
Therefore by the Sobolev embedding theorem, we have $S^{2,2}\subset S^{1,4}$
and%
\begin{equation}
||\lambda -\overline{\lambda }||_{S^{1,4}}\leq C(T).  \label{1.3}
\end{equation}

Now using Theorem \ref{T2.1} (pseudohermitian Moser inequality), we get 
\begin{equation*}
\int_{M}e^{4(\lambda -\overline{\lambda })}d\mu _{0}\leq C\exp (C||\lambda -%
\overline{\lambda }||_{S^{1,4}})\leq C(T).
\end{equation*}
Together with $\int_{M}e^{4\lambda }d\mu _{0}$ being invariant under the
flow, we conclude that%
\begin{equation}
C\geq \overline{\lambda }\geq -C(T)  \label{eqn1.9}
\end{equation}
(the upper bound is obtained by observing that $\int \lambda d\mu _{0}\leq
\int e^{4\lambda }d\mu _{0}).$ By (\ref{eqn1.8}) and (\ref{eqn1.9}), we
finally obtain 
\begin{equation*}
||\lambda ||_{S^{2k,2}}\leq C(k,T)
\end{equation*}%
for all $0\leq t<T.$

\endproof%

\textbf{Proof of Theorem \ref{Th1.1}:}

When the torsion is zero, we can write $P_{0}=\Box _{0}\overline{\Box _{0}}$ 
$=$ $\overline{\Box _{0}}\Box _{0}$ ($\Box _{0}$ and $\overline{\Box _{0}}$
commute) where $\Box _{0}$ $\equiv $ $\Delta _{0}+iT_{0}$ is the Kohn
Laplacian (acting on functions) with respect to $(J,\theta _{0}),$ and $%
T_{0} $ is an infinitesimal $CR$ diffeomorphism. This implies that ($M,J)$
is embeddable (see Theorem 2.1 in \cite{le}). Therefore $\Box _{0}$ and $%
\overline{\Box _{0}}$ have subelliptic estimates on the orthogonal
complements of $Ker(\Box _{0})$ and $Ker(\overline{\Box _{0}})$ in $L^{2},$
denoted as ($Ker(\Box _{0}))^{\perp }$ and ($Ker(\overline{\Box _{0}}%
))^{\perp },$ respectively (\cite{k1}, \cite{k2}). Since $P_{0}$ commutes
with $\Box _{0}$ and $Ker(\overline{\Box _{0}})$ $\subset $ $KerP_{0}$, we
have

\begin{equation}
\Box _{0}\varphi \in (KerP_{0})^{\perp }\subset (Ker(\overline{\Box _{0}}%
))^{\perp }  \label{eqn1.10}
\end{equation}%
for $\varphi \in (KerP_{0})^{\perp }.$ We can then estimate

\begin{eqnarray}
||\varphi ||_{S^{k+4,2}} &\leq &C||\Box _{0}\varphi ||_{S^{k+2,2}}
\label{eqn1.11} \\
&\leq &C||\overline{\Box _{0}}\Box _{0}\varphi ||_{S^{k,2}}=C||P_{0}\varphi
||_{_{S^{k,2}}}  \notag
\end{eqnarray}%
for $\varphi \in (KerP_{0})^{\perp }$ ($\subset $($Ker(\Box _{0}))^{\perp })$
by (\ref{eqn1.10}). By (\ref{eqn1.11}), we mean that $P_{0}$ is subelliptic
on $(KerP_{0})^{\perp }$. It follows that there exists a unique $C^{\infty }$
smooth solution $\lambda ^{\perp }$ of (\ref{11}) for a short time (noting
that $P_{0}$ is essentially positive by Theorem \ref{t2.1}, which explains
the negative sign in (\ref{11})). On the other hand, we can apply the
contraction mapping principle to show the short time existence of a unique $%
C^{\infty }$ smooth solution $\lambda _{\ker }$ to (\ref{11a}) (observing
that $\lambda _{\ker }$ satisfies (\ref{1.1})). The long time solution then
follows from Proposition \ref{p3.1}, the Sobolev embedding theorem for $%
S^{k,2},$ and the standard argument for extending the solution at the
maximal time $T.$

\endproof%

\section{\textbf{A Priori Estimates and Asymptotic Convergence}}

We can obtain an a priori $S^{2,2}$ estimate under some analytic assumption
replacing the torsion free condition. Also a condition on the background $Q$%
-curvature will assure that the bound is uniform, i.e., independent of the
time. Under the same condition together with the vanishing torsion, we will
then have the higher-order estimates with unform bounds. Therefore we are
able to prove the asymptotic convergence of solutions of (\ref{2}).

\begin{lemma}
\label{4.1} Suppose ($Q_{0})_{\ker }=0.$ Then the solution $\lambda ^{\perp
} $ of (\ref{11}) satisfies the following a priori estimate%
\begin{equation}
\text{ }\mathcal{E(\lambda }^{\perp }\mathcal{)}\equiv \int_{M}P_{0}\lambda
^{\perp }\cdot \lambda ^{\perp }d\mu _{0}+\int_{M}Q_{0}\lambda ^{\perp }d\mu
_{0}\leq \beta ^{2}  \label{eqn4.1}
\end{equation}%
for a constant $\beta $ independent of the time.
\end{lemma}

\proof
The condition ($Q_{0})_{\ker }=0$ implies that 
\begin{equation*}
\int_{M}Q_{0}\lambda _{\ker }d\mu _{0}=0.
\end{equation*}

It follows that 
\begin{equation*}
\begin{array}{l}
\int_{M}P_{0}\lambda \cdot \lambda d\mu _{0}+\int_{M}Q_{0}\lambda d\mu _{0}
\\ 
=\int_{M}P_{0}\lambda ^{\perp }\cdot (\lambda _{\ker }+\lambda ^{\perp
})d\mu _{0}+\int_{M}Q_{0}(\lambda _{\ker }+\lambda ^{\perp })d\mu _{0} \\ 
=\int_{M}P_{0}\lambda ^{\perp }\cdot \lambda ^{\perp }d\mu
_{0}+\int_{M}Q_{0}\lambda ^{\perp }d\mu _{0}.%
\end{array}%
\end{equation*}
Hence from Lemma \ref{l3.3}, we obtain (\ref{eqn4.1}).

\endproof%

\textbf{Proof of Theorem \ref{Th1.2}:}

The essential positivity of $P_{0}$ assures that for some constant $\Upsilon
,$ there holds%
\begin{equation}
\int_{M}P_{0}\lambda ^{\perp }\cdot \lambda ^{\perp }d\mu _{0}\geq \Upsilon
\int_{M}(\lambda ^{\perp })^{2}d\mu _{0}.  \label{eqn4.2}
\end{equation}%
From the proof of Proposition \ref{p3.1}, we have (\ref{1.2}), i.e., a bound
for the $L^{2}$ norm of $\lambda ^{\perp },$ depending on the maximal time $%
T $. If in addition we assume ($Q_{0})_{\ker }=0,$ then by (\ref{eqn4.1}), (%
\ref{eqn4.2}), and the Young inequality, we obtain 
\begin{equation}
\int_{M}(\lambda ^{\perp })^{2}d\mu _{0}\leq C(\Upsilon ,Q_{0},\theta
_{0},\beta )  \label{26a}
\end{equation}%
and 
\begin{equation}
\int_{M}P_{0}\lambda \cdot \lambda d\mu _{0}=\int_{M}P_{0}\lambda ^{\perp
}\cdot \lambda ^{\perp }d\mu _{0}\leq C(\Upsilon ,Q_{0},\theta _{0},\beta )
\label{26}
\end{equation}%
for all $t\geq 0.$

Now if the condition $(\ast )$ holds on \ $(M,J,[\theta _{0}])$, that is, 
\begin{equation*}
\int_{M}|\overset{0}{\nabla ^{2}}\lambda ^{\perp }|^{2}d\mu _{0}\leq
(2+\varepsilon _{0})\int_{M}(\Delta _{0}\lambda ^{\perp })^{2}d\mu
_{0}+C(\varepsilon _{0})\int_{M}|\overset{0}{\nabla }\lambda ^{\perp
}|^{2}d\mu _{0}.
\end{equation*}%
>From this and Lemma \ref{l3.2}, we have 
\begin{equation*}
2\int_{M}P_{0}\lambda ^{\perp }\cdot \lambda ^{\perp }d\mu _{0}\geq
(1-\varepsilon _{0})\int_{M}(\Delta _{0}\lambda ^{\perp })^{2}d\mu
_{0}-C\int_{M}|\overset{0}{\nabla }\lambda ^{\perp }|^{2}d\mu _{0}.
\end{equation*}%
By applying (\ref{26a}), integrating by parts and the Young inequality, we
compute

\begin{equation}
\begin{array}{ccl}
2\int_{M}P_{0}\lambda ^{\perp }\cdot \lambda ^{\perp }d\mu _{0} & = & 
(1-\varepsilon _{0})\int_{M}(\Delta _{0}\lambda ^{\perp })^{2}d\mu
_{0}+C\int_{M}\lambda ^{\perp }\Delta _{0}\lambda ^{\perp }d\mu _{0} \\ 
& \geq & (1-\varepsilon _{0}-\varepsilon )\int_{M}(\Delta _{0}\lambda
^{\perp })^{2}d\mu _{0}-C(\varepsilon _{0},\frac{1}{\varepsilon }),%
\end{array}
\label{28}
\end{equation}%
for small $\varepsilon >0.$ Now if we choose $\varepsilon $ small enough
such that $(1-\varepsilon _{0}-\varepsilon )=\delta >0$, then 
\begin{equation}
2\int_{M}P_{0}\lambda ^{\perp }\cdot \lambda ^{\perp }d\mu _{0}\geq \delta
\int_{M}(\Delta _{0}\lambda ^{\perp })^{2}d\mu _{0}-C(\varepsilon _{0},\frac{%
1}{\varepsilon }).  \label{30}
\end{equation}%
Finally from (\ref{26}) and (\ref{30}), there exists a positive constant $%
C=C(\delta ,\Upsilon ,Q_{0},\theta _{0},\beta )$ such that%
\begin{equation*}
\int_{M}(\Delta _{0}\lambda ^{\perp })^{2}d\mu _{0}\leq C,
\end{equation*}%
and hence%
\begin{equation}
||\lambda ^{\perp }||_{S^{2,2}}\leq C  \label{34}
\end{equation}%
for all $t\geq 0$ by (\ref{26a}) (the constant $C$ depends on the maximal
time $T$ by (\ref{1.2}) if we don't assume ($Q_{0})_{\ker }=0).$

To deal with $\lambda _{\ker },$ we consider $\lambda -\bar{\lambda}.$
Without the condition ($Q_{0})_{\ker }=0,$ we have obtained the bound for $%
||\lambda ||_{S^{2,2}},$ depending on $T,$ from the later part of the proof
of Proposition \ref{p3.1}. Let us assume $(Q_{0})_{\ker }=0$ below and see
how to get an uniform bound.

Under the condition $(Q_{0})_{\ker }=0$, the equation (\ref{11a}) reads 
\begin{equation*}
\frac{\partial \lambda _{\ker }}{\partial t}=r,
\end{equation*}%
and hence we have%
\begin{equation}
\lambda _{\ker }(t,p)=\lambda _{\ker }(0,p)+\int_{0}^{t}rdt.  \label{1.6}
\end{equation}%
It follows that

\begin{equation*}
\lambda -\overline{\lambda }=(\lambda ^{\perp }-\overline{\lambda ^{\perp }}%
)+[\lambda _{\ker }(0,p)]
\end{equation*}%
and from (\ref{34}) and (\ref{26a}) that%
\begin{equation*}
||\lambda -\overline{\lambda }||_{S^{2,2}}\leq C_{1}+||\lambda ^{\perp }-%
\overline{\lambda ^{\perp }}||_{S^{2,2}}\leq C_{2}
\end{equation*}%
for all $t\geq 0.$ Therefore by the Sobolev embedding theorem for $S^{k,p}$,
we have%
\begin{equation}
||\lambda -\overline{\lambda }||_{S^{1,4}}\leq C_{3}  \label{34a}
\end{equation}%
for all $t\geq 0$ $.$

Now again by using Theorem \ref{T2.1}, we have 
\begin{equation*}
\int_{M}e^{4(\lambda -\overline{\lambda })}d\mu _{0}\leq C\exp (C||\lambda -%
\overline{\lambda }||_{S^{1,4}})\leq C_{4}
\end{equation*}%
for all $t\geq 0.\ $Together with $\int_{M}e^{4\lambda }d\mu _{0}$ being
invariant under the flow, we conclude that%
\begin{equation}
C_{5}\geq \overline{\lambda }\geq -C_{5}.  \label{1.7}
\end{equation}%
Thus 
\begin{equation*}
||\lambda ||_{S^{2,2}}\leq C
\end{equation*}%
(and also (from Theorem \ref{T2.1})%
\begin{equation}
\int_{M}e^{\alpha \lambda }d\mu _{0}\leq C  \label{51}
\end{equation}%
for all real number $\alpha $ and all $t\geq 0).$

\endproof%

Since $\Delta _{0}$ is self adjoint, we can deduce the following result by
the Young inequality.

\begin{lemma}
\label{L2} (interpolation inequality for $\Delta _{0})$ On a closed
pseudohermitian manifold $(M,J,\theta _{0}),$ given $\varepsilon >0,$ there
exists a constant $C(\varepsilon )$ such that 
\begin{equation}
\int_{M}(\Delta _{0}\varphi )^{2}d\mu _{0}\leq \varepsilon \int_{M}(\Delta
_{0}^{2}\varphi )^{2}d\mu _{0}+C(\varepsilon )\int_{M}\varphi ^{2}d\mu _{0}
\label{eqn4.2"}
\end{equation}

for all real valued functions $\varphi \in C^{\infty }(M).$
\end{lemma}

\textbf{Proof of Theorem \ref{Th1.3}:}

We compute%
\begin{equation}
\begin{array}{l}
\frac{d}{dt}\int_{M}(\Delta _{0}^{k}\lambda ^{\perp })^{2}d\mu _{0} \\ 
=2\int_{M}(\Delta _{0}^{k}\lambda ^{\perp })(\Delta _{0}^{k}\frac{\partial
\lambda ^{\perp }}{\partial t})d\mu _{0} \\ 
=-2\int_{M}(\Delta _{0}^{k}\lambda ^{\perp })[\Delta
_{0}^{k}(Q_{0}+2P_{0}\lambda ^{\perp })]d\mu _{0}\text{ \ (by (\ref{11}))}
\\ 
=-2\int_{M}(\Delta _{0}^{k}\lambda ^{\perp })[\Delta _{0}^{k}(Q_{0})]d\mu
_{0}-4\int_{M}(\Delta _{0}^{k}\lambda ^{\perp })[\Delta
_{0}^{k}(P_{0}\lambda ^{\perp })]d\mu _{0}.%
\end{array}
\label{3.5}
\end{equation}

Next we will estimate the second term in the last line of (\ref{3.5}). From
the Bochner formula (\ref{3.1}), (\ref{1.4}), and (\ref{1.5}), we compute%
\begin{equation}
\begin{array}{l}
-4\int_{M}(\Delta _{0}^{k}\lambda ^{\perp })[\Delta _{0}^{k}(P_{0}\lambda
^{\perp })]d\mu _{0} \\ 
=-4\int_{M}(\Delta _{0}^{k}\lambda ^{\perp })[\Delta _{0}^{k}(\Delta
_{0}^{2}\lambda ^{\perp }+T_{0}^{2}\lambda ^{\perp })]d\mu _{0}+E_{1} \\ 
=-4\int_{M}(\Delta _{0}^{k+1}\lambda ^{\perp })^{2}d\mu
_{0}-4\int_{M}(\Delta _{0}^{k}\lambda ^{\perp })[\Delta
_{0}^{k}T_{0}^{2}\lambda ^{\perp }]d\mu _{0}+E_{1} \\ 
=-4\int_{M}(\Delta _{0}^{k+1}\lambda ^{\perp })^{2}d\mu
_{0}-4\int_{M}(\Delta _{0}^{k}\lambda ^{\perp })[T_{0}^{2}\Delta
_{0}^{k}\lambda ^{\perp }]d\mu _{0}+E_{1}+E_{2} \\ 
=-4\int_{M}(\Delta _{0}^{k+1}\lambda ^{\perp })^{2}d\mu
_{0}+4\int_{M}(T_{0}\Delta _{0}^{k}\lambda ^{\perp })^{2}d\mu
_{0}+E_{1}+E_{2} \\ 
=2[\int_{M}|\overset{0}{\nabla ^{2}}(\Delta _{0}^{k}\lambda ^{\perp
})|^{2}d\mu _{0}-3\int_{M}(\Delta _{0}^{k+1}\lambda ^{\perp })^{2}d\mu
_{0}]+E_{1}+E_{2}+E_{3}%
\end{array}
\label{eqn4.2a}
\end{equation}%
where 
\begin{equation*}
\begin{array}{ccl}
E_{1} & = & \int_{M}(\Delta _{0}^{k}\lambda ^{\perp })\Delta
_{0}^{k}(O_{2}(\lambda ^{\perp }))d\mu _{0}, \\ 
E_{2} & = & \int_{M}(\Delta _{0}^{k}\lambda ^{\perp })O_{2k}(\lambda ^{\perp
}))d\mu _{0}, \\ 
E_{3} & = & \int_{M}(\Delta _{0}^{k}\lambda ^{\perp })O_{2k+2}(\lambda
^{\perp })d\mu _{0},%
\end{array}%
\end{equation*}%
with $O_{j}$ being some differential operator of weight $j$ (see Appendix $A$
for the definition). Using the subelliptic estimate and the Young
inequality, we can estimate the \textquotedblright error
terms\textquotedblright\ $E_{j}^{\prime }s$ by%
\begin{equation}
\begin{array}{l}
|E| \\ 
=|E_{1}+E_{2}+E_{3}| \\ 
\leq C(k,\varepsilon _{1})\int_{M}(\Delta _{0}^{k}\lambda ^{\perp })^{2}d\mu
_{0}+\varepsilon _{1}(\int_{M}(\Delta _{0}^{k+1}\lambda ^{\perp })^{2}d\mu
_{0}+\int_{M}(\lambda ^{\perp })^{2}d\mu _{0}). \\ 
\multicolumn{1}{c}{}%
\end{array}
\label{eqn4.2b}
\end{equation}%
Note that $\Delta _{0}$ (and hence $\Delta _{0}^{k})$ maps $Ker(P_{0}),$ and
hence $Ker(P_{0})^{\bot },$ into itself by assumption. So we can apply the
condition $(\ast )$ to the first term of the last line in (\ref{eqn4.2a}).
Together with substituting (\ref{eqn4.2b}) into (\ref{eqn4.2a}), we obtain%
\begin{equation}
\begin{array}{l}
-4\int_{M}(\Delta _{0}^{k}\lambda ^{\perp })[\Delta _{0}^{k}(P_{0}\lambda
^{\perp })]d\mu _{0} \\ 
\leq -2(1-\varepsilon _{0}-\frac{\varepsilon _{1}}{2})\int_{M}(\Delta
_{0}^{k+1}\lambda ^{\perp })^{2}d\mu _{0}+C(k,\varepsilon
_{1})\int_{M}(\Delta _{0}^{k}\lambda ^{\perp })^{2}d\mu _{0}+\varepsilon
_{1}\int_{M}(\lambda ^{\perp })^{2}d\mu _{0}.%
\end{array}
\label{eqn4.2c}
\end{equation}%
Taking $\varepsilon _{1}$ small so that $1-\varepsilon _{0}-\frac{%
\varepsilon _{1}}{2}$ $>$ $0.$ From the interpolation inequality (\ref%
{eqn4.2"}) with $\varphi $ $=$ $\Delta _{0}^{k-1}\lambda ^{\perp },$ we get

\begin{equation}
\int_{M}(\Delta _{0}^{k}\lambda ^{\perp })^{2}d\mu _{0}\leq \varepsilon
\int_{M}(\Delta _{0}^{k+1}\lambda ^{\perp })^{2}d\mu _{0}+C(\varepsilon
)\int_{M}(\Delta _{0}^{k-1}\lambda ^{\perp })^{2}d\mu _{0}.  \label{eqn4.2d}
\end{equation}
By choosing $\varepsilon $ (depending on $k,\varepsilon _{1})$ small enough,
we can absorb the middle term of the second line into the first term of the
second line in (\ref{eqn4.2c}). Substituting (\ref{eqn4.2d}) into (\ref%
{eqn4.2c}) and then substituting the result into (\ref{3.5}), we obtain

\begin{equation}
\begin{array}{l}
\frac{d}{dt}\int_{M}(\Delta _{0}^{k}\lambda ^{\perp })^{2}d\mu _{0} \\ 
\leq -C\int_{M}(\Delta _{0}^{k+1}\lambda ^{\perp })^{2}d\mu
_{0}+C(k)\int_{M}(\Delta _{0}^{k-1}\lambda ^{\perp })^{2}d\mu
_{0}+C\int_{M}(\lambda ^{\perp })^{2}d\mu _{0}.%
\end{array}
\label{eqn4.2e}
\end{equation}%
In view of the subellipticity of $\Delta _{0}^{k},$ we have%
\begin{equation*}
\int_{M}(\Delta _{0}^{k}\lambda ^{\perp })^{2}d\mu _{0}\leq
C(k)(\int_{M}(\Delta _{0}^{k+1}\lambda ^{\perp })^{2}d\mu
_{0}+\int_{M}(\lambda ^{\perp })^{2}d\mu _{0}).
\end{equation*}%
Therefore we can reduce (\ref{eqn4.2e}) to

\begin{eqnarray}
&&\frac{d}{dt}\int_{M}(\Delta _{0}^{k}\lambda ^{\perp })^{2}d\mu _{0}
\label{eqn4.2f} \\
&\leq &-C(k)\int_{M}(\Delta _{0}^{k}\lambda ^{\perp })^{2}d\mu
_{0}+C(k)\int_{M}(\Delta _{0}^{k-1}\lambda ^{\perp })^{2}d\mu
_{0}+C(k)\int_{M}(\lambda ^{\perp })^{2}d\mu _{0}.  \notag
\end{eqnarray}

Now we apply induction on $k\geq 2$ to (\ref{eqn4.2e}). When $k=2,$ the last
two terms of (\ref{eqn4.2f}) are bounded due to Theorem \ref{Th1.2}. So by
Lemma \ref{ODE}, we get 
\begin{equation*}
\int_{M}(\Delta _{0}^{2}\lambda ^{\perp })^{2}d\mu _{0}\leq C(T)
\end{equation*}%
and $C(T)$ is independent of $T$ if the condition ($Q_{0})_{\ker }=0$ is
imposed. By induction hypothesis, $\int_{M}(\Delta _{0}^{k-1}\lambda ^{\perp
})^{2}d\mu _{0}$ is bounded. Also recall in the proof of Proposition \ref%
{p3.1} that the essential positivity of $P_{0}$ and the energy estimate (\ref%
{9}) imply the $L^{2}$ bound of $\lambda ^{\perp }$ (see (\ref{1.2}) and (%
\ref{26a}) if, in addition, ($Q_{0})_{\ker }=0).$ So the last two terms of (%
\ref{eqn4.2f}) are bounded. Applying Lemma \ref{ODE} again, we conclude that

\begin{equation}
\int_{M}(\Delta _{0}^{k}\lambda ^{\perp })^{2}d\mu _{0}\leq C(k,T)
\label{eqn4.2g}
\end{equation}%
where $C(k,T)$ is independent of $T$ if the condition ($Q_{0})_{\ker }=0$ is
imposed. From (\ref{1.2}) ((\ref{26a}) if, in addition, ($Q_{0})_{\ker }=0$%
), (\ref{eqn4.2g}), and the subelliptic estimate, we can then have

\begin{equation*}
||\lambda ^{\perp }||_{S^{2k,2}}\leq C(k,T)
\end{equation*}%
($C(k,T)$ is independent of $T$ if the condition ($Q_{0})_{\ker }=0$ is
imposed). Observe that the argument from (\ref{eqn1.8}) to (\ref{eqn1.9}) in
the proof of Proposition \ref{p3.1} still works without torsion free
condition. We therefore have the following estimate

\begin{equation*}
||\lambda ||_{S^{2k,2}}\leq C(k,T)
\end{equation*}%
with $C(k,T)$ being independent of $T$ if the condition ($Q_{0})_{\ker }=0$
is imposed.

\endproof%

\begin{proposition}
\label{t3.5} Let $(M,J,[\theta _{0}])$ be a CR $3$-manifold with $\overset{0}%
{A}_{11}=0$ and $(Q_{0})_{\ker }=0$. There exists a constant $C=C(\Upsilon
,Q_{0},\theta _{0},\beta )\;$such that 
\begin{equation*}
||\lambda ||_{S^{4,2}}\leq C
\end{equation*}%
for all $t\geq 0.$ In particular, there holds 
\begin{equation*}
|\lambda |\leq C
\end{equation*}%
for all $t\geq 0.$ Moreover, we have 
\begin{equation}
||\lambda ||_{S^{2k,2}}\leq C(k)  \label{eqn4.2'}
\end{equation}%
for all $t\geq 0.$
\end{proposition}

\proof
Observe that $\lambda $ has an uniform $L^{2}$ bound if $(Q_{0})_{\ker }=0$
in the proof of Proposition (\ref{p3.1}), and hence the constants $C(...,T)$%
, $C(T)$, and $C(k,T)$ can be replaced by $C$ or $C(k)$ (independent of $T$)
from (\ref{1.2}) to the end of the proof.

\endproof%

\textbf{Proof of Theorem \ref{t1}:}

By Theorem \ref{Th1.1}, we have a long time solution. Starting from (\ref{2d}%
), we compute%
\begin{eqnarray}
\frac{d^{2}}{dt^{2}}\mathcal{E}(\theta ) &=&-4\int (Q_{0}+2P_{0}\lambda
)P_{0}\frac{\partial \lambda }{\partial t}d\mu _{0}  \label{eqn4.3} \\
&=&4\int (Q_{0}^{\perp }+2P_{0}\lambda )P_{0}(Q_{0}^{\perp }+2P_{0}\lambda
)d\mu _{0}\text{ by (\ref{2})}  \notag \\
&\geq &4\Upsilon \int (Q_{0}^{\perp }+2P_{0}\lambda )^{2}  \notag
\end{eqnarray}%
by essential positivity of $P_{0}$ (which is implied by the torsion free
condition; see Theorem \ref{t2.1}). Therefore $\frac{d}{dt}\mathcal{E}%
(\theta )$ is nondecreasing, and hence%
\begin{equation}
\int_{M}e^{4\lambda }Q^{2}d\mu \text{ (}\geq 0)\text{ is nonincreasing}.
\label{eqn4.4}
\end{equation}%
By (\ref{eqn4.2'}), we can find a sequence of times $t_{j}$ such that $%
\lambda _{j}$ $\equiv $ $\lambda (\cdot ,t_{j})$ converges to $\lambda
_{\infty }$ in $C^{\infty }$ topology as $t_{j}$ $\rightarrow \infty $. On
the other hand, integrating (\ref{2d}) gives%
\begin{equation}
\mathcal{E}(\lambda _{\infty })-\mathcal{E}(\lambda _{0})=-\int_{0}^{\infty
}\int_{M}e^{4\lambda }Q^{2}d\mu dt.  \label{eqn4.5}
\end{equation}%
In view of (\ref{eqn4.4}), (\ref{eqn4.5}), we obtain

\begin{equation}
0=\lim_{t\rightarrow \infty }\int_{M}e^{4\lambda }Q^{2}d\mu
=\int_{M}e^{4\lambda _{\infty }}Q_{\infty }^{2}d\mu _{\infty }
\label{eqn4.6}
\end{equation}%
where $Q_{\infty }$ denotes the $Q$-curvature with respect to $(J,\theta
_{\infty }),$ $\theta _{\infty }$ $=$ $e^{2\lambda _{\infty }}\theta _{0},$
and $d\mu _{\infty }$ $=$ $e^{4\lambda _{\infty }}d\mu _{0}.$ It follows
that $Q_{\infty }$ $=$ $0.$

In the following, we are going to prove the smooth convergence for all time.
First we want to prove that $\lambda $ converges to $\lambda _{\infty }$ in $%
L^{2}.$ Write $\lambda _{\infty }$ $=$ $\lambda _{\infty }^{\perp }+(\lambda
_{\infty })_{\ker }.$ Observe that ($||\cdot ||_{2}$ denotes the $L^{2}$
norm with respect to the volume form $d\mu _{0}$)%
\begin{eqnarray}
||\lambda _{\infty }^{\perp }-\lambda ^{\perp }||_{2} &\leq &||\lambda
_{\infty }^{\perp }-\lambda ^{\perp }||_{S^{4,2}}  \label{eqn4.7} \\
&\leq &C||2P_{0}(\lambda _{\infty }^{\perp }-\lambda ^{\perp })||_{2}  \notag
\\
&=&C||2P_{0}\lambda ^{\perp }+Q_{0}^{\perp }||_{2}  \notag
\end{eqnarray}%
by the subellipticity of $P_{0}$ on ($KerP_{0})^{\perp }$ and $0$ $=$ $%
Q_{\infty }$ $=$ $e^{-4\lambda _{\infty }}(2P_{0}\lambda _{\infty }^{\perp
}+Q_{0}^{\perp })$ (($Q_{0})_{\ker }$ $=$ $0$ by assumption). We compute%
\begin{equation}
\begin{array}{l}
\mid \mathcal{E}(\lambda _{\infty }^{\perp })-\mathcal{E}(\lambda ^{\perp
})\mid \\ 
=\mid \int_{0}^{1}\frac{d}{ds}\mathcal{E(}\lambda ^{\perp }+s(\lambda
_{\infty }^{\perp }-\lambda ^{\perp }))ds\mid \\ 
=\mid \int_{0}^{1}\int_{M}[2P_{0}(\lambda ^{\perp }+s(\lambda _{\infty
}^{\perp }-\lambda ^{\perp }))+Q_{0}^{\perp }]\cdot (\lambda _{\infty
}^{\perp }-\lambda ^{\perp })d\mu _{0}ds\mid \\ 
\leq \int_{0}^{1}||2P_{0}(\lambda ^{\perp }+s(\lambda _{\infty }^{\perp
}-\lambda ^{\perp }))+Q_{0}^{\perp }||_{2}||\lambda _{\infty }^{\perp
}-\lambda ^{\perp }||_{2}ds \\ 
\leq C_{1}||2P_{0}\lambda ^{\perp }+Q_{0}^{\perp }||_{2}^{2}%
\end{array}
\label{eqn4.8}
\end{equation}%
by the Cauchy inequality and (\ref{eqn4.7}). Let $\vartheta $ be a number
between $0$ and $\frac{1}{2}.$ It follows from (\ref{eqn4.8}) that%
\begin{equation}
\mid \mathcal{E}(\lambda _{\infty }^{\perp })-\mathcal{E}(\lambda ^{\perp
})\mid ^{1-\vartheta }\leq C_{2}||2P_{0}\lambda ^{\perp }+Q_{0}^{\perp
}||_{2}^{2(1-\vartheta )}\leq C_{2}||2P_{0}\lambda ^{\perp }+Q_{0}^{\perp
}||_{2}  \label{eqn4.9}
\end{equation}%
for $t$ large by noting that $2(1-\vartheta )$ $>$ $1$ and $||2P_{0}\lambda
^{\perp }+Q_{0}^{\perp }||_{2}$ tends to $0$ as $t\rightarrow \infty $ by (%
\ref{eqn4.6}). Next we compute%
\begin{equation}
\begin{array}{l}
-\frac{d}{dt}(\mathcal{E}(\lambda ^{\perp })-\mathcal{E}(\lambda _{\infty
}^{\perp }))^{\vartheta } \\ 
=-\vartheta (\mathcal{E}(\lambda ^{\perp })-\mathcal{E}(\lambda _{\infty
}^{\perp }))^{\vartheta -1}\frac{d}{dt}(\mathcal{E}(\lambda ^{\perp })-%
\mathcal{E}(\lambda _{\infty }^{\perp })) \\ 
=+\vartheta (\mathcal{E}(\lambda ^{\perp })-\mathcal{E}(\lambda _{\infty
}^{\perp }))^{\vartheta -1}||2P_{0}\lambda ^{\perp }+Q_{0}^{\perp }||_{2}||%
\dot{\lambda}^{\perp }||_{2} \\ 
\geq \vartheta C_{2}^{-1}||\dot{\lambda}^{\perp }||_{2}%
\end{array}
\label{eqn4.10}
\end{equation}%
by (\ref{2d}), (\ref{eqn4.9}), and noting that $\dot{\lambda}^{\perp }$ $=$ $%
-(2P_{0}\lambda ^{\perp }+Q_{0}^{\perp })$ (see (\ref{11})) and hence the
left side is nonnegative. We learned the above trick of raising the power to 
$\vartheta $ from \cite{s}. Integrating (\ref{eqn4.10}) with respect to $t$
gives

\begin{equation}
\int_{0}^{\infty }||\dot{\lambda}^{\perp }||_{2}dt<+\infty .  \label{eqn4.11}
\end{equation}

Observing that $\frac{d}{dt}(\lambda ^{\perp }-\lambda _{\infty }^{\perp })$ 
$=$ $\dot{\lambda}^{\perp }$ and $-\frac{d}{dt}||\lambda ^{\perp }-\lambda
_{\infty }^{\perp }||_{2}^{2}$ $\leq $ $C||\dot{\lambda}^{\perp }||_{2},$ we
can then deduce 
\begin{equation}
\lim_{t\rightarrow \infty }||\lambda ^{\perp }-\lambda _{\infty }^{\perp
}||_{2}^{2}=0  \label{eqn4.11'}
\end{equation}%
by (\ref{eqn4.11}). On the other hand, we can estimate $|r(t)|$ $\leq $ $%
C||2P_{0}\lambda ^{\perp }+Q_{0}^{\perp }||_{2}$ $=$ $C||\dot{\lambda}%
^{\perp }||_{2}$ by (\ref{51}). It follows from (\ref{eqn4.11}) that

\begin{equation}
\int_{0}^{\infty }|r(t)|dt<+\infty .  \label{eqn4.12}
\end{equation}%
So in view of (\ref{1.6}) and (\ref{eqn4.12}), $\lambda _{\ker }$ converges
to $(\lambda _{\infty })_{\ker }$ $=$ ($\lambda _{\ker })_{\infty }$ as $t$ $%
\rightarrow $ $\infty $ (not just a sequence of times). Since $\lambda
_{\ker }-(\lambda _{\infty })_{\ker }$ $=$ $-\int_{t}^{\infty }r(t)dt$ is a
function of time only, we also have%
\begin{equation}
\lim_{t\rightarrow \infty }||\lambda _{\ker }-(\lambda _{\infty })_{\ker
}||_{S^{k,2}}=0  \label{eqn4.13}
\end{equation}%
for any nonnegative integer $k.$ With $\lambda ^{\perp }$ replaced by $%
\lambda ^{\perp }-\lambda _{\infty }^{\perp }$ in the argument to deduce (%
\ref{eqn1.7'}), we obtain%
\begin{equation}
\frac{d}{dt}\int_{M}(\Delta _{0}^{k}(\lambda ^{\perp }-\lambda _{\infty
}^{\perp }))^{2}d\mu _{0}\leq -4\Upsilon \int_{M}(\Delta _{0}^{k}(\lambda
^{\perp }-\lambda _{\infty }^{\perp }))^{2}d\mu _{0}+C(k)||\lambda ^{\perp
}-\lambda _{\infty }^{\perp }||_{2}.  \label{eqn4.14}
\end{equation}%
By Lemma \ref{ODE} and (\ref{eqn4.11'}), we get\ 
\begin{equation}
\lim_{t\rightarrow \infty }\int_{M}(\Delta _{0}^{k}(\lambda ^{\perp
}-\lambda _{\infty }^{\perp }))^{2}d\mu _{0}=0.  \label{eqn4.15}
\end{equation}%
Hence by the subellipticity of $\Delta _{0}^{k},$ we conclude from (\ref%
{eqn4.15}) and (\ref{eqn4.11'}) that%
\begin{equation*}
\lim_{t\rightarrow \infty }||\lambda ^{\perp }-\lambda _{\infty }^{\perp
}||_{S^{2k,2}}=0.
\end{equation*}%
Together with \ref{eqn4.13}, we have proved that $\lambda $ converges to $%
\lambda _{\infty }$ smoothly as $t$ $\rightarrow $ $\infty .$

\endproof%

\section{\textbf{Essential Positivity of the }$CR$\textbf{\ Paneitz Operator}%
}

Let $(M,J,\theta $) be a closed pseudohermitian $3$-manifold with zero
torsion. We will prove the essential positivity of the CR Paneitz operator $%
P $.

Since the torsion of ($M,J,\theta )$ is zero, the Lie derivative of the $CR$
structure $J$ with respect to $T$ (characteristic vector field) is zero.
Hence $M$ admits a smooth $CR$ action of $R$, which is transverse to the
contact bundle. This implies that $(M,J)$ is embeddable by Theorem 2.1. in 
\cite{le}. Therefore the $L^{2}$-closure of the Kohn Laplacian $\Box _{b}$
has the closed range (\cite{k1}) and hence subelliptic estimates on the
orthogonal complement of the kernel of $\Box _{b}$ in $L^{2}$ (\cite{k1}, 
\cite{k2}). We use the same symbol $\Box _{b}$ to denote the $L^{2}$-closure.

On the other hand, since $\Box _{b}$ is a closed operator, $Ker(\Box _{b})$
is closed, and we have the following direct sum decomposition 
\begin{equation}
L^{2}=Ker(\Box _{b})\oplus R(\Box _{b}),  \label{1}
\end{equation}%
where $R(\Box _{b})$ is the range of $\Box _{b}$.

Let $C^{\infty }(M,C)$ be the space of all $C^{\infty }$ smooth
complex-valued functions on $M$ and $C_{K}^{\infty }=C^{\infty }(M,C)\cap
Ker(\Box _{b})$. Then we have $C^{\infty }(M,C)=C_{K}^{\infty }\oplus
(C_{K}^{\infty })^{\bot }$. By (\ref{1}) and subelliptic estimates on the
orthogonal complement of the kernel of $\Box _{b}$ in $L^{2}$, we see that $%
(C_{K}^{\infty })^{\bot }=\left\{ \Box _{b}\varphi |\ \varphi \in C^{\infty
}(M,C)\right\} $.

Let $H:C^{\infty }(M,C)\rightarrow C_{K}^{\infty }$ be the projection.
Define the Green operator $G:C^{\infty }(M,C)\rightarrow (C_{K}^{\infty
})^{\bot }$ by 
\begin{equation*}
G(\alpha )=\omega \in (C_{K}^{\infty })^{\bot },
\end{equation*}%
where $\omega $ is the unique solution of $\Box _{b}\omega =\alpha -H(\alpha
)$. It is easy to check that $G$ is symmetric and positive on $%
(C_{K}^{\infty })^{\bot }$. Moreover, since $\Box _{b}$ has subelliptic
estimates, $G$ is compact on the orthogonal complement of $Ker(\Box _{b})$.

\begin{lemma}
\label{l2.1} Let $\eta =\sup \left\{ {\left\| G\varphi \right\| :\ \ \varphi
\in C^{\infty }(M,C),\ \left\| \varphi \right\| =1,\ \varphi \in
(C_{K}^{\infty })^{\bot }}\right\} $. Then $\frac{1}{\eta }$ is an
eigenvalue of $\Box _{b}$.
\end{lemma}

\proof
Let $\left\{ \varphi _{j}\right\} \subset (C_{K}^{\infty })^{\bot }$ be a
maximizing sequence for $\eta $, that is, $\left\| \varphi _{j}\right\| =1$
and $\left\| G\varphi _{j}\right\| \rightarrow \eta $. Since $G$ is
symmetric, we have%
\begin{equation*}
\begin{array}{ccc}
||G^{2}\varphi _{j}-\eta ^{2}\varphi _{j}||^{2} & = & \left\| G^{2}\varphi
_{j}\right\| ^{2}-2\eta ^{2}\left\langle G^{2}\varphi _{j},\varphi
_{j}\right\rangle +\eta ^{4} \\ 
& \leq & \eta ^{2}\left\| G\varphi _{j}\right\| ^{2}-2\eta ^{2}\left\|
G\varphi _{j}\right\| ^{2}+\eta ^{4}.%
\end{array}%
\end{equation*}%
This inequality means that $\left\| G^{2}\varphi _{j}-\eta ^{2}\varphi
_{j}\right\| \rightarrow 0$ as $j\rightarrow 0$.

Let $\psi _{j}=G\varphi _{j}-\eta \varphi _{j}$. Since $G$ is positive on $%
(C_{K}^{\infty })^{\bot }$, we have%
\begin{equation*}
\begin{array}{ccl}
\left\langle \psi _{j},G^{2}\varphi _{j}-\eta ^{2}\varphi _{j}\right\rangle
& = & \left\langle \psi _{j},G\psi _{j}+\eta \psi _{j}\right\rangle \\ 
& = & \left\langle \psi _{j},G\psi _{j}\right\rangle +\eta \left\| \psi
_{j}\right\| ^{2} \\ 
& \geq & \eta \left\| \psi _{j}\right\| ^{2}.%
\end{array}%
\end{equation*}%
This inequality, together with the above, implies that $\left\| G\varphi
_{j}-\eta \varphi _{j}\right\| \rightarrow 0$.

On the other hand, since $G$ is compact on the orthogonal complement of the
kernel of $\Box _{b}$, there exists a subsequence of $\left\{ \varphi
_{j}\right\} $, also denoted by $\left\{ \varphi _{j}\right\} $, such that $%
\left\{ G\varphi _{j}\right\} $ converges to a function $f\in L^{2}$. We
have, for all $\varphi \in C^{\infty }(M,C)$,%
\begin{equation*}
\begin{array}{ccl}
\left\langle (\Box _{b}-\frac{1}{\eta })f,\varphi \right\rangle & = & 
\left\langle f,(\Box _{b}-\frac{1}{\eta })\varphi \right\rangle \\ 
& = & \lim_{j\rightarrow 0}\left\langle G\varphi _{j},(\Box _{b}-\frac{1}{%
\eta })\varphi \right\rangle \\ 
& = & \lim_{j\rightarrow 0}\left\langle \varphi _{j}-\frac{1}{\eta }G\varphi
_{j},\varphi \right\rangle \\ 
& = & \lim_{j\rightarrow 0}\frac{1}{\eta }\left\langle \eta \varphi
_{j}-G\varphi _{j},\varphi \right\rangle \\ 
& = & 0.%
\end{array}%
\end{equation*}%
That is, $f\in L^{2}$ is a nontrivial weak solution of $(\Box _{b}-\frac{1}{%
\eta })f=0$. Actually $f\bot Ker\Box _{b}$, so $f$ is smooth and $\frac{1}{%
\eta }$ is an eigenvalue of $\Box _{b}$.

\endproof%

\begin{lemma}
\label{l2.2} Let $(M,J,\theta )$ be a closed pseudohermitian $3$-manifold.
Suppose that $(M,J)$ is embeddable. Then $\Box _{b}$ has discrete
eigenvalues: $0=\mu _{0}<\mu _{1}\leq \mu _{2}\leq \cdots \leq \mu
_{n}\rightarrow \infty $ and the corresponding eigenfunctions $\left\{
\varphi _{i}\right| \ i\geq 1\}$ satisfying $\Box _{b}\varphi _{i}=\mu
_{i}\varphi _{i},\ \varphi _{i}\in C^{\infty }(M,C)$ can be chosen so that $%
\left\{ \varphi _{i}\right\} $ forms an orthonormal basis of $(Ker\Box
_{b})^{\bot }$ in $L^{2}$.
\end{lemma}

\proof
Let $E(\mu )$ be the eigenspace with respect to the eigenvalue $\mu $. By
the subelliptic estimates, dim$E(\mu )<\infty $ if $\mu >0$.

For $j\geq 1$, let $\eta _{j}\equiv \sup \left\{ {\left\| G\varphi \right\|
:\ \ \varphi \in C^{\infty }(M,C),\ \left\| \varphi \right\| =1,\ \varphi
\in \left( \oplus _{i=0}^{j-1}E(\frac{1}{\eta _{i}})\right) ^{\bot }}%
\right\} $. Here we denote $E(\frac{1}{\eta _{0}})=C_{K}^{\infty }$. A
similar argument as in Lemma \ref{l2.1} shows that $\frac{1}{\eta _{j}}$ is
an eigenvalue of $\Box _{b}$. Therefore we have the discrete eigenvalues $%
0=\mu _{0}<\mu _{1}\leq \mu _{2}\leq \cdots \leq \mu _{n},\ n\rightarrow
\infty $. By the subelliptic estimates, $\left\{ \mu _{n}\right\} $ can not
have a finite limit. Thus $\mu _{n}\rightarrow \infty $ as $n\rightarrow
\infty $.

Let $0=\mu _{0}<\mu _{1}\leq \mu _{2}\leq \cdots $ be the eigenvalues of $%
\Box _{b}$, where each eigenvalue is included as many times as the dimension
of its eigenspace, with a corresponding orthonormal sequence of
eigenfunctions $\left\{ \varphi _{i}\right\} $. Let $\alpha \in C^{\infty
}(M,C)$, there exists $\beta \in (C_{K}^{\infty })^{\bot }$ such that $%
G\beta =\alpha -H(\alpha )$. It is easy to check that $\left\langle \alpha
,\varphi _{i}\right\rangle \varphi _{i}=G(\left\langle \beta ,\varphi
_{i}\right\rangle \varphi _{i})$, so we have

\begin{equation*}
\begin{array}{ccl}
\left\| \alpha -\sum_{i=1}^{n}\left\langle \alpha ,\varphi _{i}\right\rangle
\varphi _{i}-H(\alpha )\right\| & = & \left\| G\beta
-\sum_{i=1}^{n}\left\langle \alpha ,\varphi _{i}\right\rangle \varphi
_{i}\right\| \\ 
& = & \left\| G\beta -G\left( \sum_{i=1}^{n}\left\langle \beta ,\varphi
_{i}\right\rangle \varphi _{i}\right) \right\| \\ 
& = & \left\| G\left( \beta -\sum_{i=1}^{n}\left\langle \beta ,\varphi
_{i}\right\rangle \varphi _{i}\right) \right\| \\ 
& \leq & \frac{1}{\mu _{n+1}}\left\| \beta -\sum_{i=1}^{n}\left\langle \beta
,\varphi _{i}\right\rangle \varphi _{i}\right\| \\ 
& \leq & \frac{1}{\mu _{n+1}}\left\| \beta \right\| .%
\end{array}%
\end{equation*}%
Then $\mu _{n}\rightarrow \infty $ implies that 
\begin{equation*}
\lim_{n\rightarrow \infty }\left\| \alpha -\sum_{i=1}^{n}\left\langle \alpha
,\varphi _{i}\right\rangle \varphi _{i}-H(\alpha )\right\| =0.
\end{equation*}

\endproof%

\begin{theorem}
\label{t2.1}If the torsion of ($M,J,\theta )$ is zero, then the CR Paneitz
operator $P$ is essentially positive.
\end{theorem}

\proof
First, zero torsion implies that (i) the $CR$ Paneitz operator $P=\Box _{b}%
\overline{\Box _{b}}$, and (ii) the Kohn Laplacian $\Box _{b}$ and $%
\overline{\Box _{b}}$ commute, so they are diagonalized simultaneously on
the finite dimensional eigenspace of $\Box _{b}$ with respect to any nonzero
eigenvalue. Therefore, we can choose an orthonormal basis $\left\{ \varphi
_{i}\right\} $ such that each eigenfunction $\in \left\{ \varphi
_{i}\right\} $ is also an eigenfunction of $\overline{\Box _{b}}$, and
hence, of $P$. We know that the eigenvalues of $\Box _{b}$ (and hence of $%
\overline{\Box _{b}}$) are all nonnegative. Therefore by Lemma \ref{l2.2}, $%
P $ is essentially positive.

\endproof%
\bigskip

\appendix

\section{ \ }

We will give a brief introduction to pseudohermitian geometry (see \cite{l1}%
, \cite{l2} for more details). Let $M$ be a closed $3$-manifold with an
oriented contact structure $\xi $. There always exists a global contact form 
$\theta $, obtained by patching together local ones with a partition of
unity. The characteristic vector field of $\theta $ is the unique vector
field $T$ such that ${\theta }(T)=1$ and $\mathcal{L}_{T}{\theta }=0$ or $d{%
\theta }(T,{\cdot })=0$. A $CR$ structure compatible with $\xi $ is a smooth
endomorphism $J:{\xi }{\rightarrow }{\xi }$ such that $J^{2}=-identity$. A
pseudohermitian structure compatible with $\xi $ is a $CR$-structure $J$
compatible with $\xi $ together with a global contact form $\theta $.

Given a pseudohermitian structure $(J,{\theta })$, we can choose a complex
vector field $Z_{1}$, an eigenvector of $J$ with eigenvalue $i$, and a
complex $1$-form ${\theta }^{1}$ such that $\{{\theta },{\theta ^{1}},{%
\theta ^{\bar{1}}}\}$ is dual to $\{T,Z_{1},Z_{\bar{1}}\}$. It follows that $%
d{\theta }=ih_{1{\bar{1}}}{\theta ^{1}}{\wedge }{\theta ^{\bar{1}}}$ for
some nonzero real function $h_{1{\bar{1}}}$. If $h_{1{\bar{1}}}$ is
positive, we call such a pseudohermitian structure $(J,{\theta })$ positive,
and we can choose a $Z_{1}$ (hence $\theta ^{1}$) such that $h_{1{\bar{1}}%
}=1 $. That is to say%
\begin{equation*}
d\theta =i\theta ^{1}\wedge \theta ^{\overline{1}}.
\end{equation*}

We will always assume our pseudohermitian structure $(J,{\theta })$ is
positive and $h_{1{\bar{1}}}=1$ throughout the paper. The pseudohermitian
connection of $(J,{\theta })$ is the connection ${\nabla }^{{\psi }.h.}$ on $%
TM{\otimes }C$ (and extended to tensors) given by%
\begin{equation*}
{\nabla }^{{\psi }.h.}Z_{1}={\omega _{1}}^{1}{\otimes }Z_{1},{\nabla }^{{%
\psi }.h.}Z_{\bar{1}}={\omega _{\bar{1}}}^{\bar{1}}{\otimes }Z_{\bar{1}},{%
\nabla }^{{\psi }.h.}T=0
\end{equation*}%
\noindent in which the $1$-form ${\omega _{1}}^{1}$ is uniquely determined
by the following equation with a normalization condition:%
\begin{equation}
\begin{array}{ccl}
d{\theta ^{1}} & = & {\theta ^{1}}{\wedge }{\omega _{1}}^{1}+{A^{1}}_{\bar{1}%
}{\theta }{\wedge }{\theta ^{\bar{1}}} \\ 
0 & = & {\omega _{1}}^{1}+{\omega _{\bar{1}}}^{\bar{1}}.%
\end{array}
\label{2.2}
\end{equation}

The coefficient ${A^{1}}_{\bar{1}}$ in (\ref{2.2}) is called the
(pseudohermitian) torsion. Since $h_{1{\bar{1}}}=1$, $A_{{\bar{1}}{\bar{1}}%
}=h_{1{\bar{1}}}{A^{1}}_{\bar{1}}={A^{1}}_{\bar{1}}$. And $A_{11}$ is just
the complex conjugate of $A_{{\bar{1}}{\bar{1}}}$. Differentiating ${\omega
_{1}}^{1}$ gives%
\begin{equation*}
d{\omega _{1}}^{1}=W{\theta ^{1}}{\wedge }{\theta ^{\bar{1}}}+2iIm(A_{11,{%
\bar{1}}}{\theta ^{1}}{\wedge }{\theta })
\end{equation*}%
\noindent where $W$ is the Tanaka-Webster curvature.

We can define the covariant differentiations with respect to the
pseudohermitian connection. For instance, $f_{1}=Z_{1}f$, $f_{1{\bar{1}}}=Z_{%
\bar{1}}Z_{1}f-{\omega _{1}}^{1}(Z_{\bar{1}})Z_{1}f$ for a (smooth) function 
$f$. We define the sub-gradient operator $\nabla _{b}$ and the sub-Laplacian
operator $\Delta _{b}$ by%
\begin{equation*}
\begin{array}{ccl}
{\nabla _{b}}f & = & f_{\overline{1}}Z_{1}+f_{1}Z_{\bar{1}}, \\ 
{\Delta _{b}}f & = & f_{1\overline{1}}+f_{{\bar{1}}1}.%
\end{array}%
\end{equation*}%
\noindent respectively. We also define the Levi form $<,>_{J,{\theta }}$ by%
\begin{equation*}
<V,U>_{J,{\theta }}=2d{\theta }(V,JU)=v_{1}u_{\bar{1}}+v_{\bar{1}}u_{1},
\end{equation*}%
\noindent for $V=v_{1}Z_{\bar{1}}+v_{\bar{1}}Z_{1}$,$U=u_{1}Z_{\bar{1}}+u_{%
\bar{1}}Z_{1}$ in $\xi $ and%
\begin{equation*}
(V,U)_{J,\theta }=\int_{M}<V,U>_{J,\theta }\theta \wedge d\theta {.}
\end{equation*}%
For a vector $X$ $\in $ $\xi ,$ we define $|X|^{2}\equiv <X,X>_{J,{\theta }%
}. $ It follows that $|\nabla _{b}f|^{2}$ $=$ $2f_{1}f_{\bar{1}}$ for a real
valued smooth function $f.$ Also the square modulus of the sub-Hessian $%
\nabla _{b}^{2}f$ of $f$ reads $|\nabla _{b}^{2}f|^{2}$ $=$ $2f_{11}f_{\bar{1%
}\bar{1}}$ $+$ $2f_{1\bar{1}}f_{\bar{1}1}.$ We recall below what the
Folland-Stein space $S^{k,p}$ is$.$ Let $D$ denote a differential operator
acting on functions. We say $D$ has weight $m,$ denoted $w(D)=m,$ if $m$ is
the smallest integer such that $D$ can be locally expressed as a polynomial
of degree $m$ in vector fields tangent to the contact bundle $\xi .$ We
define the Folland-Stein space $S^{k,p}$ of functions on $M$ by

\begin{equation*}
S^{k,p}=\{f\in L^{p}:Df\in L^{p}\text{ whenever }w(D)\leq k\}.
\end{equation*}

\noindent We define the $L^{p}$ norm of $\nabla _{b}f,$ $\nabla _{b}^{2}f$,
... to be ($\int |\nabla _{b}f|^{p}\theta \wedge d\theta )^{1/p},$ ($\int
|\nabla _{b}^{2}f|^{p}\theta \wedge d\theta )^{1/p},$ $...,$ respectively,
as usual. So it is natural to define the $S^{k,p}$ norm of $f\in S^{k,p}$ as
follows:%
\begin{equation*}
||f||_{S^{k,p}}\equiv (\sum_{0\leq j\leq k}||\nabla
_{b}^{j}f||_{L^{p}}^{p})^{1/p}.
\end{equation*}

\noindent The function space $S^{k,p}$ with the above norm is a Banach space
for $k\geq 0,$ $1<p<\infty .$ There are also embedding theorems of Sobolev
type. For instance, $S^{2,2}\subset S^{1,4}$ (for $\dim M$ $=$ $3$). We
refer the reader to, for instance, \cite{fs2} and \cite{fo} for more
discussions on these spaces.

\bigskip


\bigskip

\begin{thebibliography}{CC1}
\bibitem[A]{a} T. Aubin, \textit{Some Nonlinear Problems in Riemannian
Geometry}, Springer-Verlag, Berlin, 1998.

\bibitem[B]{b} S. Brendle, Global Existence and Convergence for a
Higher-Order Flow in Conformal Geometry, Ann. Math., 158 (2003), 323-343.

\bibitem[C]{c} P. T. Chru\'{s}ciel, Semi-Global Existence and Convergence of
solutions of the Robinson-Trautman (2-Dimensional Calabi) Equation, Commun.
Math. Phys. 137 (1991), 289-313.

\bibitem[C1]{c1} S.- C. Chang, The $Q$-Curvature Flow on a Closed $3$%
-Manifold of Positive $Q$-Curvature, preprint.

\bibitem[C2]{c2} S.- C. Chang, Recent Developments on the Calabi Flow,
Contemporary Mathematics, 367, 17-42, A.M.S., 2005.

\bibitem[C3]{c3} S.- C. Chang, Global Existence and Convergence of Solutions
of the Calabi Flow on Riemann Surfaces of Genus $g\geq 2$, J. Math. Kyoto
Univ. Vol. 40, No. 2 (2000), 363-377.

\bibitem[C4]{c4} S.- C. Chang, The $2$-dimensional Calabi Flow, Nagoya Math.
J., to appear.

\bibitem[CC1]{cc1} S.- C. Chang and J.- H. Cheng, The Harnack Estimate for
the Yamabe Flow on $CR$ Manifolds of Dimension $3$, Annals of Global
Analysis and Geometry, 21 (2002), 111-121.

\bibitem[CH]{ch} S.- S. Chern and R. Hamilton, On Riemannian Metrics Adapted
to Three-dimensional Contact Manifolds, Lecture Notes in Math., 1111,
279-305, Springer-Verlag, 1984.

\bibitem[Chi]{chi} H.-L. Chiu, The Sharp Lower Bound for the First Positive
Eigenvalue of a Sub-Laplacian on a Three-Dimensional Pseudo-hermitian
Manifold, submitted.

\bibitem[CL]{cl} W. S. Cohn and G. Lu, Best Constants for Moser-Trudinger
Inequalities on the Heisenberg Group, Indiana Univ. Math. J. 50 (2001),
1567-1591.

\bibitem[CLe]{cle} J.- H. Cheng and J. M. Lee, The Burns-Epstein Invariant
and Deformation of $CR$ Structures, Duke Math. J., 60 (1990), 221-254.

\bibitem[CS]{cs} S.- C. Chen and M.- C. Shaw, \textit{Partial Differential
Equations in several comlex variables}, Studies in Advan. Math., 19, AMS/IP,
2001.

\bibitem[CW]{cw} S.- C. Chang and C.- T. Wu, The Fourth-Order $Q$-Curvature
Flow on Closed $3$-Manifolds, Nagoya Math. J., to appear.

\bibitem[Fa]{fa} F. Farris, An Intrinsic Construction of Fefferman's $CR$
Metric, Pacific J. Math., 123 (1986), 33-45.

\bibitem[FH]{fh} C. Fefferman and K. Hirachi, Ambient Metric Construction of 
$Q$-Curvature in Conformal and $CR$ Geometries, to appear in Math. Res.
Lett..

\bibitem[Fo]{fo} G. B. Folland, Subelliptic Estimates and Function Spaces on
Nilpotent Lie Groups, Arkiv for Mat. 13 (1975), 161-207.

\bibitem[FS1]{fs1} G. B. Folland and E. M. Stein, \textit{Hardy spaces on
homogeneous groups}, Princeton U. Press\textit{, }1980.

\bibitem[FS2]{fs2} G. B. Folland and E. M. Stein, Estimates for the $\bar{%
\partial}_{b}$ Complex and Analysis on the Heisenberg Group, Comm. Pure
Appl. Math., 27 (1974), 429-522.

\bibitem[GG]{gg} A. R. Gover and C. R. Graham, $CR$ Invariant Powers of the
Sub-Laplacian, preprint.

\bibitem[GL]{gl} C. R. Graham and J. M. Lee, Smooth Solutions of Degenerate
Laplacians on Strictly Pseudoconvex Domains, Duke Math. J., 57 (1988),
697-720.

\bibitem[H]{h} K. Hirachi, Scalar Pseudo-hermitian Invariants and the Szeg%
\"{o} Kernel on $3$-dimensional $CR$ Manifolds, Lecture Notes in Pure and
Appl. Math. 143, pp. 67-76, Dekker, 1992.

\bibitem[JL]{jl} D. Jerison and J. M. Lee, The Yamabe Problem on $CR$
Manifolds, J. Diff. Geom. 25 (1987), 167-197.

\bibitem[K1]{k1} J. J. Kohn, Estimates for $\overline{\partial }_{b}$ on
Compact Pseudoconvex $CR$ Manifolds, Proc. of Symposia in Pure Math., 43
(1985), 207--217.

\bibitem[K2]{k2} J. J. Kohn, The Range of the Tangential Cauchy-Riemann
Operator, Duke. Math. J., 53 (1986) 525--545.

\bibitem[L1]{l1} J. M. Lee, Pseudo-Einstein Structure on $CR$ Manifolds,
Amer. J. Math. 110 (1988), 157-178.

\bibitem[L2]{l2} J. M. Lee, The Fefferman Metric and Pseudohermitian
Invariants, Trans. A.M.S. 296 (1986), 411-429.

\bibitem[Le]{le} L. Lempert, On Three-Dimensional Cauchy-Riemann Manifolds.
J. of Amer. Math. Soc.,\textit{\ }\textbf{5} (1992), 923--969.

\bibitem[P]{p} S. Paneitz, A Quartic Conformally Covariant Differential
Operator for Arbitrary Pseudo-Riemannian Manifolds, preprint, 1983.

\bibitem[S]{s} L. Simon, Asymptotics for a Class of Nonlinear Evolution
Equations, with Applications to Geometric Problems, Ann. of Math., 118
(1983), 525-571.

\bibitem[SC]{sc} S\'{a}nchez-Calle A., Fundamental Solutions and Geometry of
the Sum of Squares of Vector Fields, Invent. Math., 78 (1984), 143--160.
\end{thebibliography}
\end{document}